\documentclass[10pt]{article}

\usepackage{amsthm}
\usepackage{amssymb}
\usepackage{mathrsfs}
\usepackage{empheq}
\usepackage{dsfont}
\usepackage[applemac]{inputenc}
\usepackage{hyperref}

\newcommand{\1}{\mathds{1}}

\newcommand{\N}{\mathbb{N}}

\newcommand{\R}{\mathbb{R}}
\newcommand{\C}{\mathbb{C}}

\newcommand{\vH}{{\cal H}}

\newcommand{\Dr}{\mathscr{D}}

\newcommand{\Lr}{\mathscr{L}}

\newcommand{\vphi}{\varphi}
\newcommand{\eps}{\varepsilon}

\newcommand{\dsp}{\displaystyle}
\newcommand{\ovl}{\overline}
\newcommand{\udl}{\underline}
\newcommand{\vlim}{\lim\limits}

\newcommand{\vsup}{\sup\limits}

\newcommand{\vint}{\int\limits}
\newcommand{\vsum}{\sum\limits}
\newcommand{\inj}{\hookrightarrow}
\newcommand{\tends}{\longrightarrow}
\newcommand{\weak}{\rightharpoonup}
\newcommand{\wt}{\widetilde}

\newcommand{\loc}{\mathrm{loc}}

\renewcommand{\b}{\mathrm{b}}
\newcommand{\co}{\mathrm{c}}
\renewcommand{\d}{\mathrm{d}}

\newcommand{\GN}{\mathrm{GN}}
\newcommand{\vi}{\mathrm{i}}
\newcommand{\w}{{\textsl w}}
\renewcommand{\le}{\leqslant}
\renewcommand{\ge}{\geqslant}
\renewcommand{\Re}{\mathrm{Re}}
\renewcommand{\Im}{\mathrm{Im}}
\newcommand{\bs}{\boldsymbol}
\newcommand{\p}{\prime}

\newcommand{\eqdef}{\stackrel{\mathrm{def}}{=}}
\DeclareMathOperator{\supp}{supp}

\numberwithin{equation}{section}

\newtheorem{thm}{Theorem}[section]
\newtheorem{prop}[thm]{Proposition}
\newtheorem{cor}[thm]{Corollary}
\newtheorem{lem}[thm]{Lemma}

\theoremstyle{definition}
\newtheorem{rmk}[thm]{Remark}
\newtheorem{defi}[thm]{Definition}

\newtheorem{assum}[thm]{Assumption}

\newenvironment{proof*}{\noindent{\bf Proof.}}{\qed}
\newenvironment{vproof}[1]{\noindent{\bf Proof #1}}{\qed}

\title{\huge \sc Finite time extinction for a critically damped Schrödinger equation with a sublinear nonlinearity}
\author{\sc Pascal Bégout$^*$ and Jes\'us Ildefonso D{\'{\i}}az$^\dagger$}
\date{}

\begin{document}

\maketitle

\begin{gather*}
\begin{array}{cc}
                         ^*\mbox{Institut de Mathématiques de Toulouse}	&	\;^\dagger\mbox{Instituto de Matem\'atica Interdisciplinar}		\\
                                          \mbox{Université Toulouse I Capitole}	&	\mbox{Universidad Complutense de Madrid}				\\
                                             \mbox{1, Esplanade de l’Université}	&	\mbox{Plaza de las Ciencias, 3}						\\
                                  \mbox{31080 Toulouse Cedex 6, FRANCE}	&	\mbox{28040 Madrid, SPAIN}							\\
\bigskip \\
\mbox{
{\footnotesize E-mail\:: \href{mailto:Pascal.Begout@math.cnrs.fr}{\udl{\texttt{Pascal.Begout@math.cnrs.fr}}}}}
&
\mbox{
{\footnotesize E-mail\:: \href{mailto:jidiaz@ucm.es}{\udl{\texttt{jidiaz@ucm.es}}}}
}
\end{array}
\end{gather*}

\begin{abstract}
This paper completes some previous studies by several authors on the finite time extinction for nonlinear Schr\"{o}dinger equation when the nonlinear damping term corresponds to the limit cases of some ``saturating non-Kerr law'' $F(|u|^2)u=\frac{a}{\eps+(|u|^2)^\alpha}u,$ with $a\in\C,$ $\eps\ge0,$ $2\alpha=(1-m)$ and $m\in[0,1).$ Here we consider the sublinear case $0<m<1$ with a critical damped coefficient: $a\in\C$ is assumed to be in the set $D(m)=\big\{z\in\C; \; \Im(z)>0 \text{ and } 2\sqrt m\Im(z)=(1-m)\Re(z)\big\}.$ Among other things, we know that this damping coefficient is critical, for instance, in order to obtain the monotonicity of the associated operator (see the paper by Liskevich and Perel$^\p$muter~\cite{MR1224619} and the more recent study by Cialdea and Maz$^\p$ya~\cite{MR2155899}). The finite time extinction of solutions is proved by a suitable energy method after obtaining appropiate a priori estimates. Most of the results apply to non-necessarily bounded spatial domains.
\end{abstract}

{\let\thefootnote\relax}
{\let\thefootnote\relax}

\tableofcontents

\baselineskip .6cm

\section{Introduction}
\label{introduction}

In this paper, we are interested in the existence, uniqueness and finite time extinction of solutions of the damped nonlinear Schr\"{o}dinger equation
\begin{empheq}[left=\empheqlbrace]{align}
	\label{nls}
	\vi\frac{\partial u}{\partial t}+\Delta u+V(x)u+a|u|^{-(1-m)}u=f(t,x),	&	\text{ in } (0,\infty)\times\Omega,				\\
	\label{nlsb}
	u_{|\partial\Omega}=0,									&	\text{ on } (0,\infty)\times\partial\Omega,	\dfrac{}{}	\\
	\label{u0}
	u(0)= u_0,												&	\text{ in } \Omega,
\end{empheq}
where $\vi^2=-1,$ $0<m<1,$ $a\in\C$ satisfies
\begin{gather*}
2\sqrt m\,\Im(a)=(1-m)\Re(a)>0,
\end{gather*}
$\Omega\subseteq\R^N$ non-necessarily bounded, $f\in L^1_\loc\big([0,\infty);L^2(\Omega)\big),$ $V\in L^1_\loc(\Omega;\R)$ and $u_0\in L^2(\Omega).$ The finite time extinction of the solutions was first establihed in Carles and Gallo~\cite{MR2765425} in the following case: $a=\vi,$ $0\le m<1,$ $V=0,$ $f=0$ and $\Omega$ is a compact manifold without boundary. In the same paper, existence and uniquess of $H^1$ and $H^2$-solutions, in the sense quite close to the Definitions~\ref{defsol} and \ref{defsolH2} below, are shown by using a compactness method. In Carles and Ozawa~\cite{MR3306821}, the authors obtain existence and uniqueness of $H^1$ and $H^2$-solutions for some additional nonlinearities. The closest case to our study is the following: $a=\vi\lambda,$ $0\le m\le1,$ $V=-\vsum_{j=1}^N\omega_j|x_j|^2,$ $\lambda,\omega_1,\ldots\omega_N>0,$ $f=0,$ $\Omega=\R^N$ and $N\in\{1,2\}$ with also $\frac12\le m\le1,$ if $N=2.$
\medskip
\\
In this paper, we are interested by establishing existence and uniqueness results for the equation~\eqref{nls} with $m\in(0,1),$ set in an arbitrary open subset $\Omega\subseteq\R^N$ and for the largest range of $a$ as possible. For $m\in[0,1],$ let us introduce the following sets of complex numbers:
\begin{align}
\label{Cm}
& C(m)=\Big\{z\in\C; \; \Im(z)>0 \text{ and } 2\sqrt m\Im(z)\ge(1-m)|\Re(z)|\Big\},		\\
\label{Dm}
& D(m)=\Big\{z\in\C; \; \Im(z)>0 \text{ and } 2\sqrt m\Im(z)=(1-m)\Re(z)\Big\}.
\end{align}
Note that $D(0)=C(0),$ $D(1)=\emptyset$ and
\begin{align*}
& C(0)=\Big\{z\in\C; \; \Re(z)=0 \text{ and } \Im(z)>0\Big\},	\\
& C(1)=\Big\{z\in\C; \; \Im(z)>0\Big\}.
\end{align*}
Here and after, for $z\in\C,$ $\Re(z),$ $\Im(z)$ and $\ovl z$ denote the real part, the imaginary part and the conjugate of $z,$ respectively. Existence and uniqueness have been established in the following cases.

1) For $0<m<1.$

\indent\indent
a) $a\in C(m),$ $V=0$ and $|\Omega|<\infty$ (\cite{MR4053613});

\indent\indent
b) $a\in C(m)\setminus D(m),$ $V=0$ and $\Omega=\R^N$ (\cite{MR4098330});

\indent\indent
c) $a\in C(m)\setminus D(m)$ (\cite{MR4340780}).

2) For $m\in\{0,1\}.$

\indent\indent
a) $m=0,$ $a\in C(0)$ and $|\Omega|<\infty$ (\cite{MR4340780});

\indent\indent
b) $m=1,$ $a\in C(1)$ and $V=0$ (\cite{MR4053613});

\indent\indent
c) $m=1$ and $a\in C(1)$ (\cite{MR4340780}).
\\
In a nutshell, the cases
\begin{gather*}
\Omega \text{ arbitrary, } 0<m<1 \text{ and } a\in C(m)\setminus D(m),	\\
\Omega \text{ arbitrary, } m=1 \text{ and } a\in C(1),					\\
|\Omega|<\infty, \;  m=0 \text{ and } a\in C(0),
\end{gather*}
have been completely treated. It remains the cases
\begin{gather}
\label{r1}
\Omega \text{ arbitrary, } 0<m<1 \text{ and } a\in D(m),				\\
\label{r2}
|\Omega|=\infty, \;  m=0 \text{ and } a\in C(0),
\end{gather}
where \eqref{r1} can be viewed as a limit case:
\begin{gather*}
\text{for } 0<m<1 \text{ and } a\in D(m), \; a=\vlim_{\left\{\substack{\wt a\to a \hfill \\ \wt a\in C(m)\setminus D(m)}\right.}\wt a.
\end{gather*}
In this paper, we are interested by~\eqref{r1}, while~\eqref{r2} could be the subject of a future work.
\medskip
\\
A fundamental argument in our approach is the fact that
\begin{gather*}
\eqref{nls} \iff \frac{\d u}{\d t}+Au=f.
\end{gather*}
Then, we are interested in the application of the abstract theory of maximal monotone operator to the corresponding operator on the Hilbert space $L^2(\Omega).$
\medskip
\\
In \cite{MR4053613} it was directly shown that $(D(A),A)$ is maximal monotone by using the embedding $L^p(\Omega)\inj L^2(\Omega),$ for any $p>2,$ once we assume $|\Omega|<\infty.$
\medskip
\\
A different point of view was followed in \cite{MR4098330}. It was shown there that $(D(A),A)$ is maximal monotone in the following way. First, constructing solutions compactly supported in $H^2(\R^N)$ to $(A+I)u=F$ with help of the results in \cite{MR2876246,MR3315701}. Second, obtaining a priori estimates in $H^2$ with \cite[Lemma~4.2]{MR4098330}. Third, showing that $(D(A),A)$ is maximal monotone by approximations with solutions compactly supported.
\medskip
\\
A second different argument was used in \cite{MR4340780}. First, approximating $(D(A),A)$ by a nice maximal monotone operator $(D(A_\eps),A_\eps).$ Second, obtaining a priori estimates in $H^2$ with \cite[Lemma~4.2]{MR4098330}. Third, passing to the limit in the equation $(I+A_\eps)u_\eps=F,$ to prove that $(D(A),A)$ is maximal monotone.
\medskip
\\
It is important to point out that if $a\in D(m)$ then \cite[Lemma~4.2]{MR4098330} is no more valid. Then, a third argument could be apply by approximating $(D(A),A)$ by a nice maximal monotone operator $(D(A_\eps),A_\eps)$ and, by passing to the limit, to show that  $(D(A),A)$ is maximal monotone in antoher way than in~\cite{MR4340780}, by choosing $D(A)$ bigger than that of~\cite{MR4340780} (see \ref{scr3} in Section~\ref{scr}).
\medskip
\\
Notice that we are interested in the case in which $\Re(a)>0$. When $a\in\R$ and $a>0$ the general nonlinear Schr\"odinger equation is called as ``the focusing case'' (see, e.g. the exposition made by Weinstein in~\cite{MR3409633}, p.41-79) then, depending of the value of the power in the nonlinearity, global
existence in time or blow up in finite time occur. Here, by the contrary, $a\in\C$ and $\Im(a)\neq0$. As a consequence, the conservations laws (mass and energy) are broken and then the solution, which is global in time, goes to $0$ at infinity in the $L^2$-norm (the so called mass of the solution). 
\medskip
\\
This paper is organized as follows. In Section~\ref{exiuni} we present several results on the existence and uniqueness of different types of solutions. The statements of the results on finite time extinction and asymptotic behaviour of solutions are collected in Section~\ref{finite}. The proofs of the existence of solutions theorems are given in Section~\ref{proofexi}. The special case of $H^2$-solutions is considered in Section~\ref{H2sol}. Section~\ref{proofext} contains the proofs of the finite time extinction and asymptotic behavior theorems. Finally, some open problems and other remarks are collected in Section~\ref{scr}.
\medskip
\\
To end this introduction, we collect here some notations which will be used along with this paper. Let $\Omega$ be an open subset of $\R^N.$ Unless if specified, all functions are complex-valued $(H^1(\Omega)\eqdef H^1(\Omega;\C),$ etc) and all the vector spaces are considered over the field $\R.$ For $p\in[1,\infty],$ $p^\prime$ is the conjugate of $p$ defined by $\frac{1}{p}+\frac{1}{p^\prime}=1.$ For a (real)  Banach space $X,$ we denote by $X^\star\eqdef\Lr(X;\R)$ its topological dual and by $\langle\: . \; , \: . \:\rangle_{X^\star,X}\in\R$ the $X^\star-X$ duality product. When $X$ (respectively, $X^\star)$ is endowed of the weak topology $\sigma(X,X^\star)$ (respectively, the weak$\star$ topology $\sigma(X^\star,X)),$ it is denoted by $X_\w$ (respectively, by $X_{\w\star}).$ For $p\in(0,\infty],$ $u\in L^p_\loc\big([0,\infty);X\big)$ means that $u\in L^p_\loc\big((0,\infty);X\big)$ and for any $T>0,$ $u_{|(0,T)}\in L^p\big((0,T);X\big).$ In the same way, we will use the notation $u\in W^{1,p}_\loc\big([0,\infty);X\big).$ The scalar product in $L^2(\Omega)$ between two functions $u,v$ is, $(u,v)_{L^2(\Omega)}=\Re\int_{\Omega}u(x)\ovl{v(x)}\d x.$ $L^0(\Omega)$ is the space of measurable functions $u:\Omega\tends\C$ such that $|u|<\infty,$ almost eveywhere in $\Omega.$ Auxiliary positive constants will be denoted by $C$ and may change from a line to another one. Also for positive parameters $a_1,\ldots,a_n,$ we shall write $C(a_1,\ldots,a_n)$ to indicate that the constant $C$ depends only and continuously on $a_1,\ldots,a_n.$

\section{Existence and uniqueness of solutions}
\label{exiuni}

The following assumptions will be needed to construct solutions.
\begin{assum}
\label{ass1}
We assume the following.
\begin{gather}
\label{O}
\Omega \text{ is any nonempty open subset of } \R^N,	\\
\label{m}
0<m<1,										\\
\label{a}
a\in D(m),										\\
\label{V}
V\in L^\infty(\Omega;\R)+L^{p_V}(\Omega;\R),
\end{gather}
where,
\begin{gather}
\label{pV}
p_V=
\begin{cases}
2,							&	\text{if } N=1,	\\
2+\beta, \text{ for some } \beta>0,	&	\text{if } N=2,	\\
N,							&	\text{if } N\ge3.
\end{cases}
\end{gather}
\medskip
\end{assum}

\begin{rmk}
\label{rmka}
The assumption~\eqref{pV} on $p_V$ is needed to have that $Vu\in L^2(\Omega),$ for any $u\in H^1_0(\Omega)$ (see~\eqref{lemVL2} below). The proof relies on Hölder's inequality and the Sobolev embeddings (see \cite[Lemma~4.1]{MR4340780} for the complete proof). But the same proof works if $V$ satisfies the assumption
\begin{gather}
\label{Vq}
V\in L^\infty(\Omega;\R)+L^{q_V}(\Omega;\R),
\end{gather}
where
\begin{gather}
\label{qV}
q_V\in
\begin{cases}
[2,\infty],	&	\text{if } N=1,	\\
(2,\infty], 	&	\text{if } N=2,	\\
[N,\infty],	&	\text{if } N\ge3,
\end{cases}
\end{gather}
which seems to be weaker since if $V$ satisfies~\eqref{V}--\eqref{pV} then it satisfies~\eqref{Vq}--\eqref{qV} with $q_V=p_V.$ But actually, it is not. Indeed, we claim that,
\begin{gather*}
L^\infty(\Omega;\R)+L^{q_V}(\Omega;\R)\subset L^\infty(\Omega;\R)+L^{p_V}(\Omega;\R),
\end{gather*}
where it is understood that $p_V=q_V,$ if $N=2$ and $q_V<\infty.$ The claim beeing clear if $q_V=\infty,$ we are brought back to the case where $N\neq2$ and $q_V<\infty.$ Let then $V=V_1+V_2\in L^\infty(\Omega;\R)+L^{q_V}(\Omega;\R),$ where $q_V$ satisfies~\eqref{qV}. To prove the claim, it is sufficient to show that $V_2\in L^\infty(\Omega;\R)+L^{p_V}(\Omega;\R).$ Since $p_V\le q_V,$ we have that,
\begin{gather*}
\left|V_2\1_{\{|V_2|>1\}}\right|\le|V_2|^\frac{q_V}{p_V}\in L^{p_V}(\Omega;\R),
\end{gather*}
so that,
\begin{gather*}
|V_2|=\left|V_2\1_{\{|V_2|\le1\}}\right|+\left|V_2\1_{\{|V_2|>1\}}\right|\in L^\infty(\Omega;\R)+L^{p_V}(\Omega;\R).
\end{gather*}
Hence the claim.
\end{rmk}

\noindent
Here and after, we shall always identify $L^2(\Omega)$ with its topological dual. Let us recall some important results of functional analysis. Let $E$ and $F$ be locally convex Hausdorff topological vector spaces. If $E\overset{e}{\inj}F$ with dense embedding then $F^\star\overset{e^\star}{\inj}E^\star,$ where $e^\star$ is the transpose of $e:$
\begin{gather}
\label{dualtranspose}
\forall L\in F^\star, \; \forall x\in E, \; \langle e^\star(L),x\rangle_{E^\star,E}=\langle L,e(x)\rangle_{F^\star,F}.
\end{gather}
If, furthermore, $E$ is reflexive then the embedding $F^\star\overset{e^\star}{\inj}E^\star$ is dense. In most of the cases, $e$ is the identity function, so that $e^\star$ is nothing else but the restriction to $E$ of continuous linear forms on $F.$ In particular, if $X$ is a Banach space such that  $X\inj L^p(\Omega)$ with dense embedding, for some $p\in[1,\infty),$ then $L^{p^\p}(\Omega)\inj X^\star$ and for any $u\in L^{p^\p}(\Omega)$ and $v\in X,$
\begin{gather}
\label{dualg}
\langle u,v\rangle_{X^\star,X}=\langle u,v\rangle_{L^{p^\p}(\Omega),L^p(\Omega)}
=\Re\vint_{\Omega} u(x)\ovl{v(x)}\d x.
\end{gather}
For more details, see Trèves~\cite[Corollary~5, p.188; Corollary, p.199; Theorem~18.1, p.184]{MR2296978} and~\cite{MR4521439}.
Let $A_1$ and $A_2$ be two Banach spaces such that $A_1,A_2\subset\vH$ for some Hausdorff topological vector space $\vH.$ Then $A_1\cap A_2$ and $A_1+A_2$ are Banach spaces where,
\begin{gather*}
\|a\|_{A_1\cap A_2}=\max\big\{\|a\|_{A_1},\|a\|_{A_2}\big\}	\; \text{ and } \;
\|a\|_{A_1+A_2}=\inf_{\left\{\substack{a=a_1+a_2 \hfill \\ (a_1,a_2)\in A_1\times A_2}\right.}\Big(\|a_1\|_{A_1}+\|a_2\|_{A_2}\Big).
\end{gather*}
If, in addition, $A_1\cap A_2$ is dense in both $A_1$ and $A_2$ then,
\begin{gather}
\label{dual}
\big(A_1\cap A_2\big)^\star=A_1^\star+A_2^\star \; \text{ and } \; \big(A_1+A_2\big)^\star=A_1^\star\cap A_2^\star.
\end{gather}
See, for instance, Bergh and Löfström~\cite{MR0482275} (Lemma~2.3.1 and Theorem~2.7.1). Let $1<q<\infty$ and $X$ be a Banach space such that  $X\inj L^2(\Omega)$ with dense embedding. We have by~Bégout and D\'{\i}az~\cite[Lemma~A.4]{MR4053613},
\begin{gather}
\label{Xcon}
L^q_\loc\big([0,\infty);X\big)\cap W^{1,q^\p}_\loc\big([0,\infty);X^\star\big)\inj C\big([0,\infty);L^2(\Omega)\big).
\end{gather}
Let $Y$ be a Banach space such that $\Dr(\Omega)\inj Y$ with dense embedding. Then,
\begin{gather}
\label{injDp}
L^1_\loc\big((0,\infty);Y^\star\big)\inj\Dr^\p\big((0,\infty)\times\Omega\big).
\end{gather}
See, for instance, Droniou~\cite[Lemme~2.6.1]{droniou}. Finally, another result which will be useful is the following (Strauss~\cite[Theorem~2.1]{MR0205121}). Let $X\inj\Dr^\p(\Omega)$ be a reflexive Banach space. Let $I$ be an interval and $u\in C\big(\ovl I;\Dr^\p(\Omega)\big).$ If $u\in L^\infty(I;X)$ then,
\begin{gather}
\label{rmkweakcon}
\forall t\in\ovl I, \; u(t)\in X \; \text{ and } \; u\in C_\w\big(\ovl I;X\big).
\end{gather}
Here and after, $C_\w(\ovl I;X)$ denotes the space of (weakly) continuous functions from $\ovl I$ to $X_\w.$
\medskip \\
We recall the definition of solution (\cite{MR4098330,MR4053613}).

\begin{defi}
\label{defsol}
Assume~\eqref{O}, \eqref{V} and \eqref{pV}. Let $a\in\C,$  $0<m\le1,$ $f\in L^1_\loc\big([0,\infty);L^2(\Omega)\big)$ and $u_0\in L^2(\Omega).$ We shall say that $u$ is an $H^1_0$-\textit{solution of} \eqref{nls}--\eqref{u0}, if $u$ satisfies the following assertions.
\begin{enumerate}
\item
\label{defsol1}
We have,
\begin{gather*}
u\in L^{m+1}_\loc\big([0,\infty);X\big)\cap W^{1,\frac{m+1}m}_\loc\big([0,\infty);X^\star\big)
\inj C\big([0,\infty);L^2(\Omega)\big),
\end{gather*}
with $X=H^1_0(\Omega)\cap L^{m+1}(\Omega).$
\item
\label{defsol3}
$u$ satisfies~\eqref{nls} in $\Dr^\p\big((0,\infty)\times\Omega\big).$
\item
\label{defsol4}
$u(0)=u_0,$ in $L^2(\Omega).$
\end{enumerate}
We shall say that $u$ is an $L^2$-\textit{solution} or a \textit{weak solution} of \eqref{nls}--\eqref{u0} is there exists,
\begin{gather}
\label{fn}
(f_n,u_n)_{n\in\N}\subset L^1_\loc\big([0,\infty);L^2(\Omega)\big)\times C\big([0,\infty);L^2(\Omega)\big),
\end{gather}
such that for any $n\in\N,$ $u_n$ is an $\bs{H^1_0}$\textbf{-solution} of \eqref{nls}--\eqref{nlsb} where the right-hand side member of \eqref{nls} is $f_n,$ and if
\begin{gather}
\label{cv}
(f_n,u_n)\xrightarrow[n\tends\infty]{L^1((0,T);L^2(\Omega))\times C([0,T];L^2(\Omega))}(f,u),
\end{gather}
for any $T>0.$
\end{defi}

\begin{rmk}
\label{rmkdefsol}
Let us comment the Definition~\ref{defsol}.
\begin{enumerate}
\item
\label{rmkdefsol1}
In~\cite{MR4098330,MR4053613,MR4340780}, there is also a notion of $H^2$-solutions. Such solutions $u$ satisfy Properties~\ref{defsol1}--\ref{defsol4} of Definition~\ref{defsol} with, additionally, $u\in W^{1,\frac{m+1}m}_\loc\big([0,\infty);L^2(\Omega)+L^\frac{m+1}m(\Omega)\big)$ and $\Delta u(t)\in L^2(\Omega),$ for almost every $t>0$ (\cite[Definition~4.1]{MR4053613}). Unfortunately, we are not able to construct such solutions because of the lack of a priori estimates of solutions in the $H^2$-norm. Indeed, theses estimates are obtained by a rotation of $a\in C(m)\setminus D(m)$ in the complex plane, to get $a\longmapsto\wt a\in C(m).$ The crucial tool is Lemma~4.2 in Bégout~\cite{MR4098330}, which is no more valid if $a\in D(m)$ (read the proof of Bégout~\cite[Corollary~4.5]{MR4098330} to see how this lemma is applied). As a consequence, we had to modify the notion of $L^2$-solutions. Indeed, in our paper, an $L^2$-solution is a limit of $H^1_0$-solutions while in~\cite{MR4098330,MR4053613,MR4340780}, it is a limit of $H^2$-solutions. Despite this definition which seems to be weakened, such solutions do not lose any property. Indeed, the conditions \eqref{fn} and \eqref{cv} to be an $L^2$-solution are common to these four papers. As a consequence, we have not changed the terminology here. Finally, notice that $H^2$-solutions exist in the special case in which $\Omega$ has a finite measure (see Theorem~\ref{thmstrongH2} below).
\item
\label{rmkdefsol2}
The boundary condition $u(t)_{|\partial\Omega}=0$ is implicitely included in the assumption $u(t)\in H^1_0(\Omega),$ for the $H^1_0$-solutions. For the $L^2$-solutions, this has to be understood in a generalized sense by using the limit of $H^1_0$-solutions.
\end{enumerate}
\medskip
\end{rmk}

\noindent
We give an improved result from the previous paper~\cite{MR4340780} on how weak solutions satisfy~\eqref{nls} and recall a continuous dependence result.

\begin{prop}
\label{propsolL2}
Assume~\eqref{O}, \eqref{V} and \eqref{pV}. Let $0<m\le1,$ $a\in\C$ and $f\in L^1_\loc\big([0,\infty);L^2(\Omega)\big).$ Let $u$ be a weak solution to \eqref{nls}. Let $(f_n,u_n)_{n\in\N}$ satisfy~\eqref{cv}, where each $u_n$ is an $H^1_0$-solution to~\eqref{nls}--\eqref{nlsb} with $f_n$ instead of $f.$ Then,
\begin{gather}
\label{propsolL21}
u\in W^{1,1}_\loc\big([0,\infty);H^{-2}(\Omega)+L^\frac2m(\Omega)\big),
\end{gather}
and $u$ solves~\eqref{nls} in $L^1_\loc\big([0,\infty);H^{-2}(\Omega)+L^\frac2m(\Omega)\big)$ and so in $\Dr^\p\big((0,\infty)\times\Omega\big).$ In adddition,
\begin{gather}
\label{propsolL22}
u_n\xrightarrow[n\to\infty]{W^{1,1}((0,T);H^{-2}(\Omega)+L^\frac2m(\Omega))}u.
\end{gather}
for any $T>0.$
\end{prop}

\begin{prop}[\textbf{Uniqueness and continuous dependance}]
\label{propdep}
Let Assumption~$\ref{ass1}$ be fulfilled, let $f,\wt f\in L^1_\loc\big([0,\infty);L^2(\Omega)\big)$ and $X=H^1_0(\Omega)\cap L^{m+1}(\Omega).$ Finally, let
\begin{gather}
\label{propdephyp}
u,\wt u\in L^p_\loc\big([0,\infty);X\big)\cap W^{1,p^\p}_\loc\big([0,\infty);X^\star\big)\inj C\big([0,\infty);L^2(\Omega)\big),
\end{gather}
for some $1<p<\infty,$ be solutions in $\Dr^\p\big((0,\infty)\times\Omega\big)$ to,
\begin{gather*}
\vi u_t+\Delta u+Vu+a|u|^{-(1-m)}u=f,		\\
\vi\wt {u_t}+\Delta\wt u+V\wt u+a|\wt u|^{-(1-m)}\wt u=\wt f ,
\end{gather*}
respectively. Then,
\begin{gather}
\label{estthmweak}
\|u(t)-\wt u(t)\|_{L^2(\Omega)}\le\|u(s)-\wt u(s)\|_{L^2(\Omega)}+\vint_s^t\|f(\sigma)-\wt f(\sigma)\|_{L^2(\Omega)}\d\sigma,
\end{gather}
for any $t\ge s\ge0.$ Finally,~\eqref{estthmweak} also holds true for the weak solutions.
\end{prop}

\begin{thm}[\textbf{Existence and uniqueness of $\bs{L^2}$-solutions}]
\label{thmweak}
Let Assumption~$\ref{ass1}$ be fulfilled and let $f\in L^1_\loc\big([0,\infty);L^2(\Omega)\big).$ Then for any $u_0\in L^2(\Omega),$ there exists a unique weak solution $u$ to \eqref{nls}--\eqref{u0}. In addition,
\begin{gather}
\label{Lm}
u\in L^{m+1}_\loc\big([0,\infty);L^{m+1}(\Omega)\big),	\\
\label{L2+}
\dfrac12\|u(t)\|_{L^2(\Omega)}^2+\Im(a)\dsp\vint_s^t\|u(\sigma)\|_{L^{m+1}(\Omega)}^{m+1}\d\sigma
			\le\dfrac12\|u(s)\|_{L^2(\Omega)}^2+\Im\dsp\iint\limits_{s\;\Omega}^{\text{}\;\;t}f(\sigma,x)\,\ovl{u(\sigma,x)}\,\d x\,\d\sigma,
\end{gather}
for any $t\ge s\ge0.$ If $|\Omega|<\infty$ then the inequality in~\eqref{L2+} is an equality.
\end{thm}

\begin{rmk}
\label{rmkthmweak}
Using~\eqref{estthmweak}--\eqref{L2+} and Hölder's inequality, uniform continuous dependance with respect to the initial data and the right hand side member of~\eqref{nls} may be obtain in
\begin{gather*}
C_\b\big([0,\infty);L^2(\Omega)\big)\cap L^\frac{p(1-m)}{2-p}\big((0,\infty);L^p(\Omega)\big),
\end{gather*}
for any $p\in(m+1,2).$ See~\cite[Remark~2.5]{MR4098330} for more details.
\end{rmk}

\begin{thm}[\textbf{Additional regularity in $\bs{H^1_0}$ for weak solutions}]
\label{thmstrongaH1}
Let Assumption~$\ref{ass1}$ be fulfilled with additionally  $V\in W^{1,\infty}(\Omega;\R)+W^{1,p_V}(\Omega;\R).$ Let $f\in L^1_\loc\big([0,\infty);H^1_0(\Omega)\big).$ Then for any $u_0\in H^1_0(\Omega),$ the weak solution $u$ satisfies, additionally, that
\begin{gather}
\begin{cases}
\label{thmstrongaH11}
u\in C\big([0,\infty);L^2(\Omega)\big)\cap C_\w\big([0,\infty);H^1_0(\Omega)\big),	\medskip \\
u\in W^{1,1}_\loc\big([0,\infty);H^{-1}(\Omega)+L^\frac{m+1}m(\Omega)\big),
\end{cases}
\end{gather}
and $u$ satisfies \eqref{nls} in $ L^1_\loc\big([0,\infty),H^{-1}(\Omega)+L^\frac{m+1}m(\Omega)\big).$ Furthermore, $u$ verifies,
\begin{gather}
\label{thmstrongaH12}
\|u(t)\|_{H^1_0(\Omega)}\le\left(\|u(s)\|_{H^1_0(\Omega)}+\vint_s^t\|f(\sigma)\|_{H^1_0(\Omega)}\d\sigma\right)e^{C\|\nabla V\|_{L^\infty+L^{p_V}}(t-s)},
\end{gather}
for any $t\ge s\ge0,$ where $C=C(N,\beta).$ Finally, if $\nabla V=0$ then $u$ satisfies the better estimate below.
\begin{gather}
\label{thmstrongaH13}
\|\nabla u(t)\|_{L^2(\Omega)}\le\|\nabla u(s)\|_{L^2(\Omega)}+\vint_s^t\|\nabla f(\sigma)\|_{L^2(\Omega)}\d\sigma,
\end{gather}
for any $t\ge s\ge0.$
\end{thm}

\noindent
If $u$ is a weak solution given by Theorem~\ref{thmstrongaH1} and if, in addition, $f\in L^\frac{m+1}m_\loc\big((0,\infty);X^\star\big),$ where $X=H^1_0(\Omega)\cap L^{m+1}(\Omega),$ then $u$ becomes an $H^1_0$-solution, as shows the following result.

\begin{thm}[\textbf{Existence and uniqueness of $\bs{H^1_0}$-solutions -- I}]
\label{thmstrongH1}
Let Assumption~$\ref{ass1}$ be fulfilled with additionally  $V\in W^{1,\infty}(\Omega;\R)+W^{1,p_V}(\Omega;\R).$ Let
\begin{gather}
\label{fH1}
f\in L^1_\loc\big([0,\infty);H^1_0(\Omega)\big)\cap L^\frac{m+1}m_\loc\big([0,\infty);H^{-1}(\Omega)+L^\frac{m+1}m(\Omega)\big).
\end{gather}
Then for any $u_0\in H^1_0(\Omega),$ there exists a unique $H^1_0$-solution $u$ to \eqref{nls}--\eqref{u0}. Furthermore, the map $t\longmapsto\|u(t)\|_{L^2(\Omega)}^2$ belongs to $W^{1,1}_\loc\big([0,\infty);\R\big)$ and we have,
\begin{gather}
\label{L2}
\frac12\frac{\d}{\d t}\|u(t)\|_{L^2(\Omega)}^2+\Im(a)\|u(t)\|_{L^{m+1}(\Omega)}^{m+1}=\Im\vint_{\Omega}f(t,x)\,\ovl{u(t,x)}\,\d x,
\end{gather}
for almost every $t>0.$
\end{thm}

\begin{thm}[\textbf{Existence and uniqueness of $\bs{H^1_0}$-solutions -- II}]
\label{thmstrongaH2}
Let Assumption~$\ref{ass1}$ be fulfilled and $f\in W^{1,1}_\loc\big([0,\infty);L^2(\Omega)\big).$ Then for any
\begin{gather*}
u_0\in H^1_0(\Omega)\cap L^{m+1}(\Omega) \text{ for which } \Delta u_0+a|u_0|^{-(1-m)}u_0\in L^2(\Omega),
\end{gather*}
there exists a unique $H^1_0$-solution $u$ to \eqref{nls}--\eqref{u0}. Furthermore,
\begin{gather*}
u \text{ satisfies }\eqref{nls} \text { in } L^\infty_\loc\big([0,\infty);H^{-1}(\Omega)+L^\frac{m+1}m(\Omega)\big)
\end{gather*}
as well as the following properties.
\begin{enumerate}
\item
\label{thmstrongaH21}
$u\in C_\w\big([0,\infty);H^1_0(\Omega)\cap L^{m+1}(\Omega)\big)\cap W^{1,\infty}_\loc\big([0,\infty);L^2(\Omega)\big).$
\item
\label{thmstrongaH22}
For any $t\ge s\ge0,$
\begin{empheq}[left=\empheqlbrace]{align}
\label{strongaH21}
	&	\|u(t)-u(s)\|_{L^2(\Omega)}\le\|u_t\|_{L^\infty((s,t);L^2(\Omega))}|t-s|,		\frac{}{} \\
\label{strongaH22}
	&	\|u(t)\|_{L^2(\Omega)}\le A(t),									\frac{}{} \\
\label{strongaH23}
	&	\|u_t\|_{L^\infty((0,t);L^2(\Omega))}\le B(t),					\frac{}{} \\
\label{strongaH24}
	&	\|\nabla u(t)\|_{L^2(\Omega)}^2+\Im(a)\|u(t)\|_{L^{m+1}(\Omega)}^{m+1}\le C(t)A(t),
\end{empheq}
where,
\begin{align*}
&	A(t)=\|u_0\|_{L^2(\Omega)}+\int_0^t\|f(s)\|_{L^2(\Omega)}\d s,									\\
&	B(t)=\|\Delta u_0+Vu_0+ag(u_0)-f(0)\|_{L^2(\Omega)}+\int_0^t\|f^\p(\sigma)\|_{L^2(\Omega)}\d\sigma,	\\
&	C(t)=C\left(A(t),B(t),\|f(t)\|_{L^2(\Omega)},\|V_1\|_{L^\infty(\Omega)},\|V_2\|_{L^{p_V}(\Omega)},N,m,\beta\right).
\end{align*}
\item
\label{thmstrongaH23}
The map $t\longmapsto\|u(t)\|_{L^2(\Omega)}^2$ belongs to $W^{1,\infty}_\loc\big([0,\infty);\R\big)$  and~\eqref{L2} holds for almost every $t>0.$
\item
\label{thmstrongaH24}
If $f\in W^{1,1}\big((0,\infty);L^2(\Omega)\big)$ then $u\in L^\infty\big((0,\infty);H^1_0(\Omega)\cap L^{m+1}(\Omega)\big)\cap W^{1,\infty}\big((0,\infty);L^2(\Omega)\big).$
\end{enumerate}
\end{thm}

\begin{rmk}
\label{rmkf0}
Below are some comments about Theorem~\ref{thmstrongaH2}.
\begin{enumerate}
\item
\label{rmkf01}
The solution $u$ obtained in Theorem~\ref{thmstrongaH2} could be called an \textit{almost} $H^2$\textit{-solution} since it verifies all the conditions of Definition~\ref{defsolH2} below, except the property,
\begin{gather}
\label{rmkf011}
\text{for almost every } t>0, \; \Delta u(t)\in L^2(\Omega),
\end{gather}
which need not be satisfied (\cite[Definition~4.1]{MR4053613}). It merely satisfies that,
\begin{gather*}
\text{for almost every } t>0, \; \Delta u(t)\in L^2_\loc(\Omega).
\end{gather*}
The property~\eqref{rmkf011} may be obtained in the particular case in which $\Omega$ has a finite measure (see Theorem~\ref{thmstrongH2} below).
\item
\label{rmkf02}
Since $f\in W^{1,1}_\loc\big([0,\infty);L^2(\Omega)\big)\inj C\big([0,\infty);L^2(\Omega)\big),$ $f(0)$ in the function $B$ makes sense.
\item
\label{rmkf03}
For any $p\in\left(m+1,\frac{2N}{N-2}\right)$ $(p\in(m+1,\infty]$ if $N=1),$
\begin{gather*}
u\in C^{0,\alpha}\big([0,\infty);L^p(\Omega)\big)	\;
\Big(u\in C^{0,\alpha}_\b\big([0,\infty);L^p(\Omega)\big), \text{ if } f\in W^{1,1}\big((0,\infty);L^2(\Omega)\big)\Big),
\end{gather*}
where $\alpha=\frac{2N-p(N-2)}{2p}$ if $p\ge2,$ and $\alpha=2\frac{p-(1+m)}{p(1-m)}$ if $p\le2.$ Indeed, this comes from Property~\ref{thmstrongaH21} and \eqref{strongaH21}, with also Gagliardo-Nirenberg's inequality, if $p>2,$ and Hölder's inequality, if $p<2.$
\end{enumerate}
\end{rmk}

\section{Finite time extinction and asymptotic behavior}
\label{finite}

\begin{assum}
\label{assN1}
Assumption~\ref{ass1} holds true and $u_0\in H^1_0(\Omega).$ We have that $\big(f\in L^1\big((0,\infty);H^1_0(\Omega)\big)$ and $\nabla V=0\big)$ or $\big(f\in W^{1,1}\big((0,\infty);L^2(\Omega)\big),$ $u_0\in L^{m+1}(\Omega)$ with $\Delta u_0+a|u_0|^{-(1-m)}u_0\in L^2(\Omega)\big),$ and $u$ is the unique solution to \eqref{nls}--\eqref{u0} given by Theorems~\ref{thmweak} or \ref{thmstrongaH2}. Finally, there exists a $T_0\in[0,\infty)$ such that
\begin{gather}
\label{f}
\text{for almost every } t>T_0, \; f(t)=0.
\end{gather}
\end{assum}

\section*{Asymptotic behavior of the $\bs{L^2}$-solutions}

\begin{thm}
\label{thm0w}
Let Assumption~$\ref{ass1}$ be fulfilled, $f\in L^1\big((0,\infty);L^2(\Omega)\big),$ $u_0\in L^2(\Omega)$ and let $u$ be the unique weak solution to \eqref{nls}--\eqref{u0} given by Theorem~$\ref{thmweak}.$ Then,
\begin{gather*}
\vlim_{t\nearrow\infty}\|u(t)\|_{L^2(\Omega)}=0.
\end{gather*}
\end{thm}

\section*{Finite time extinction and asymptotic behavior of the $\bs{H^1_0}$-solutions}

\begin{thm}[\textbf{Finite time extinction and time decay estimates}]
\label{thmextN1}
Let Assumption~$\ref{assN1}$ be fulfilled.
\begin{enumerate}
\item
\label{thmextN11}
If $N=1$ then
\begin{gather}
\label{01}
\forall t\ge T_\star, \; \|u(t)\|_{L^2(\Omega)}=0,
\end{gather}
for some,
\begin{gather}
\label{T*N1}
T_0\le T_\star\le C\|u(T_0)\|_{L^2(\Omega)}^\frac{1-m}2\|\nabla u\|_{L^\infty((0,\infty);L^2(\Omega))}^\frac{1-m}2+T_0,
\end{gather}
for some $C=C(\Im(a),m).$
\item
\label{thmextN12}
If $N=2$ then for any $t\ge T_0,$
\begin{gather}
\label{thmrtdH11}
\|u(t)\|_{L^2(\Omega)}\le\|u(T_0)\|_{L^2(\Omega)}e^{-C(t-T_0)},
\end{gather}
where $C=C(\|\nabla u\|_{L^\infty((0,\infty);L^2(\Omega))},\Im(a),m).$
\item
\label{thmextN13}
If $N\ge3$ then for any $t\ge T_0,$
\begin{gather}
\label{thmrtdH12}
\|u(t)\|_{L^2(\Omega)}\le\dfrac{\|u(T_0)\|_{L^2(\Omega)}}
{\left(1+C\|u(T_0)\|_{L^2(\Omega)}^\frac{(1-m)(N-2)}2(t-T_0)\right)^\frac2{(1-m)(N-2)}},
\end{gather}
where $C=C(\|\nabla u\|_{L^\infty((0,\infty);L^2(\Omega))},\Im(a),N,m).$
\item
\label{thmextN14}
If $N=1,$ $f\in  L^1\big((0,\infty);H^1_0(\Omega)\big)\cap L^\frac{m+1}m\big((0,\infty);H^{-1}(\Omega)+L^\frac{m+1}m(\Omega)\big)$ and $\nabla V=0$ then there exists $\eps_\star=\eps_\star(|a|,m)$ satisfying the following property. If
\begin{gather}
\label{thmextN1e1}
\begin{cases}
\|u_0\|_{L^2(\Omega)}^{2(1-\delta_1)}\le\eps_\star T_0,							\medskip \\
\|\nabla u_0\|_{L^2(\Omega)}+\|\nabla f\|_{L^1((0,\infty);L^2(\Omega))}\le\eps_\star,		\medskip \\
\|f(t)\|_{L^2(\Omega)}^2\le\eps_\star\big(T_0-t\big)_+^\frac{2\delta_1-1}{1-\delta_1},
\end{cases}
\end{gather}
for almost every $t>0,$ where $\delta_1=\frac{3+m}4,$ then~\eqref{01} holds true with $T_\star=T_0.$
\end{enumerate}
\end{thm}

\section{Proofs of the existence of solutions}
\label{proofexi}

Before proving the results of Section~\ref{exiuni}, we recall some results of our previous paper we will need. Here and in the rest of the paper, we shall use the following notations and conventions. Unless if specified, we assume~\eqref{O}--\eqref{m}. Since $\big||z|^{-(1-m)}z\big|=|z|^m,$ we extend by continuity at $z=0$ the map $z\longmapsto|z|^{-(1-m)}z$ by setting,
\begin{gather*}
|z|^{-(1-m)}z=0, \text{ if } z=0.
\end{gather*}
Let $\eps\ge0.$ For any $u\in L^0(\Omega)$ and almost every $x\in\Omega,$ we define
\begin{align}
\label{ge}
& g_\eps^m(u)(x)=(|u(x)|^2+\eps)^{-\frac{1-m}2}u(x),	\; 0\le m\le1, \\
\label{gg}
& g(u)(x)=g_0^m(u)(x).
\end{align}
Let $p\in[1,\infty).$ We have that for any $u,v\in L^p(\Omega),$
\begin{gather}
\label{g}
\|g_0^m(u)-g_0^m(v)\|_{L^\frac{p}m(\Omega)}\le3\|u-v\|_{L^p(\Omega)}^m,
\end{gather}
In particular, $g_0^m\in C\big(L^p(\Omega);L^\frac{p}m(\Omega)\big)$ and $g_0^m$ is bounded on bounded sets. Finally, if $\eps>0$ then $g_\eps^m\in C\big(L^2(\Omega);L^2(\Omega)\big)$ and $g_\eps^m$ is bounded on bounded sets. See~\cite[Lemma~4.3]{MR4340780}.
\\
Now, let us define the operator $(A_\eps^m,D(A_\eps^m))$ on $L^2(\Omega)$ by,
\begin{gather*}
\begin{cases}
D(A_\eps^m)=\left\{u\in H^1_0(\Omega); \Delta u\in L^2(\Omega)\right\}, \medskip \\
A_\eps^mu=-\vi\Delta u-\vi Vu-\vi ag_\eps^m(u), \; \forall u\in D(A_\eps^m).
\end{cases}
\end{gather*}

\noindent
We recall the following result.

\begin{lem}[\textbf{\cite[Corollary~5.11]{MR4340780}}]
\label{lemAe}
Assume \eqref{O}. Let $0\le m<1$ and $a\in C(m).$ Then for any $\eps>0,$ $(A_\eps^m,D(A_\eps^m))$ is maximal monotone on $L^2(\Omega)$ with dense domain.
\end{lem}

\noindent
Let $V=V_1+V_2\in L^\infty(\Omega;\R)+L^{p_V}(\Omega;\R),$ where $p_V$ is given by~\eqref{pV}. Then for any $u\in L^2(\Omega),$ $Vu\in H^{-1}(\Omega)$ and for any $u\in H^1_0(\Omega),$ $Vu\in L^2(\Omega).$ There exists $C=C(N,\beta)$ such that the following holds. Let $u\in H^1_0(\Omega)$ and $v\in L^2(\Omega).$ We have,
\begin{align}
\label{lemVH-1}
&	\|Vv\|_{H^{-1}(\Omega)}\le C\|V\|_{L^\infty(\Omega)+L^{p_V}(\Omega)}\|v\|_{L^2(\Omega)},	\\
\label{lemVL2}
&	\|Vu\|_{L^2(\Omega)}\le C\|V\|_{L^\infty(\Omega)+L^{p_V}(\Omega)}\|u\|_{H^1_0(\Omega)},	\\
\label{lemdualV}
&	\langle Vv,u\rangle_{H^{-1}(\Omega),H^1_0(\Omega)}=(v,Vu)_{L^2(\Omega)},				\\
\label{lemV1}
&	\|V_1u\|_{L^2(\Omega)}\le\|V_1\|_{L^\infty(\Omega)}\|u\|_{L^2(\Omega)},					\\
\label{lemV2}
&	\|V_2u\|_{L^2(\Omega)}\le C\rho^{1-\gamma}\|V_2\|_{L^{p_V}(\Omega)}^{2-\gamma}\|u\|_{L^2(\Omega)}^\gamma+\frac1\rho\|\nabla u\|_{L^2(\Omega)}^2,
\end{align}
for any $\rho>0,$ where $\gamma=\gamma(N,\beta)\in[0,1).$ See \cite[Lemmas~4.1 and 4.2]{MR4340780}.
\\
Let us recall that for any $u\in H^1_0(\Omega)$ such that $\Delta u\in L^2(\Omega),$ we have
\begin{gather}
\label{nablau}
\|\nabla u\|_{L^2(\Omega)}^2\le\|\Delta u\|_{L^2(\Omega)}\|u\|_{L^2(\Omega)}.
\end{gather}
Finally, to prove Theorem~\ref{thmstrongaH2}, we introduce the following operator $(A,D(A))$ on $L^2(\Omega).$
\begin{gather}
\label{D}
\begin{cases}
D(A)=\left\{u\in H^1_0(\Omega)\cap L^{m+1}(\Omega); \Delta u+ag(u)\in L^2(\Omega)\right\}, \medskip \\
Au=-\vi\Delta u-\vi Vu-\vi ag(u), \, \forall u\in D(A).
\end{cases}
\end{gather}

\noindent
We have the following.

\begin{lem}
\label{lemA}
Assume~\eqref{O}--\eqref{a}. The operator $(A,D(A))$ is maximal monotone on $L^2(\Omega)$ with dense domain.
\end{lem}

\noindent
Before proving Lemma~\ref{lemA}, we give three results we will need. Lemma~\ref{lemAmax} below is stated in a more general case (in terms of $m$ and $a)$ because its proof is totally unchanged and we think that it will be of interest for a future work.

\begin{lem}
\label{lemAmax}
Assume \eqref{O}. Let $0\le m<1$ and $a\in C(m).$ Let $F\in L^2(\Omega).$ Then there exist $u\in H^1_0(\Omega)\cap L^{m+1}(\Omega),$ with $Vu\in L^2(\Omega),$ and $u_{\eps_n}\in D(A_{\eps_n}^m)$ $(n\in\N),$ where $(\eps_n)_{n\in\N}\subset(0,\infty)$ is a decreasing sequence converging toward $0,$ satisfying the following properties. For each $n\in\N,$ $u_{\eps_n}$ is the unique solution to,
\begin{gather}
\label{lemAmax1}
-\vi\Delta u_{\eps_n}-\vi Vu_{\eps_n}-\vi ag_{\eps_n}^m(u_{\eps_n})+u_{\eps_n}=F, \text{ in } L^2(\Omega).
\end{gather}
Furthermore, we have that,
\begin{gather}
\label{lemAmax2}
\vsup_{n\in\N}\|u_{\eps_n}\|_{H^1_0(\Omega)}+\sup_{n\in\N}\|Vu_{\eps_n}\|_{L^2(\Omega)}<\infty,					\\
\label{lemAmax3}
\Im(a)\vint_\Omega(|u_{\eps_n}|^2+{\eps_n})^{-\frac{1-m}2}|u_{\eps_n}|^2\d x+\|u_{\eps_n}\|_{L^2(\Omega)}^2\le\|F\|_{L^2(\Omega)}^2,
\end{gather}
for any $n\in\N.$ Finally,
\begin{align}
\label{lemAmax4}
&	u_{\eps_n}\xrightarrow[n\to\infty]{\Dr^\p(\Omega)}u,		\\
\label{lemAmax5}
&	Vu_{\eps_n}\xrightarrow[n\to\infty]{\Dr^\p(\Omega)}Vu,	\\
\label{lemAmax6}
&	u_{\eps_n}\xrightarrow[n\to\infty]{\text{a.e.\,in }\Omega}u.
\end{align}
\end{lem}

\begin{proof*}
Let $F\in L^2(\Omega).$ Let $\eps>0.$ By Lemma~\ref{lemAe} and Brezis~\cite[Proposition~2.2]{MR0348562}, there exists a unique solution $u_\eps\in D(A_\eps^m)$ to \eqref{lemAmax1} satisfying $\|u_\eps\|_{L^2(\Omega)}\le\|F\|_{L^2(\Omega)}.$ We take the $L^2$-scalar product of~\eqref{lemAmax1} with $u_\eps$ and then with $\vi u_\eps.$ We get that,
\begin{gather}
\label{demlemAmax1}
\Im(a)\vint_\Omega(|u_\eps|^2+\eps)^{-\frac{1-m}2}|u_\eps|^2\d x+\|u_\eps\|_{L^2(\Omega)}^2=\Re\vint_\Omega F\,\ovl{u_\eps}\d x,		\\
\label{demlemAmax2}
\|\nabla u_\eps\|_{L^2(\Omega)}^2-\vint_\Omega V|u_\eps|^2\d x-\Re(a)\vint_\Omega(|u_\eps|^2+\eps)^{-\frac{1-m}2}|u_\eps|^2\d x
=\Im\vint_\Omega F\,\ovl{u_\eps}\d x.
\end{gather}
Applying Cauchy-Schwarz's inequality to \eqref{demlemAmax1}, we obtain \eqref{lemAmax3}, for any sequence $\eps_n\searrow0.$ We multiply \eqref{demlemAmax1} by $\frac{\Re(a)_+}{\Im(a)},$ we sum the result with \eqref{demlemAmax2} and we still apply Cauchy-Schwarz's inequality. It follows that,
\begin{gather}
\label{demlemAmax3}
\|\nabla u_\eps\|_{L^2(\Omega)}^2\le\left(1+\frac{\Re(a)_+}{\Im(a)}\right)\|F\|_{L^2(\Omega)}^2+\|Vu_\eps\|_{L^2(\Omega)}\|F\|_{L^2(\Omega)}.
\end{gather}
By~\eqref{lemV1} and \eqref{lemV2}, there exists $C=C(N,\beta)$ such that,
\begin{gather}
\label{demlemAmax4}
\|Vu_\eps\|_{L^2(\Omega)}\le C\left(\|V_1\|_{L^\infty(\Omega)}+\|V_2\|_{L^{p_V}(\Omega)}^{2-\gamma}\right)\|F\|_{L^2(\Omega)}+\frac1{2\|F\|_{L^2(\Omega)}}\|\nabla u_\eps\|_{L^2(\Omega)}^2.
\end{gather}
With help of~\eqref{lemVL2}, \eqref{lemAmax3}, \eqref{demlemAmax3} and \eqref{demlemAmax4}, we obtain~\eqref{lemAmax2}, also for any sequence $\eps_n\searrow0.$ As a consequence, there exist $u\in H^1_0(\Omega)$ and a decreasing sequence $(\eps_n)_{n\in\N}\subset(0,\infty)$ converging toward $0$ such that, by~\eqref{lemdualV}, $Vu\in L^2(\Omega)$ and such that \eqref{lemAmax4}--\eqref{lemAmax5} hold true. By the compact embedding $H^1_0(\Omega)\inj L^2_\loc(\Omega)$ and the diagonal procedure, up to a subsequence, we get \eqref{lemAmax6}. Finally, it follows from \eqref{lemAmax3}, \eqref{lemAmax6} and Fatou's Lemma that $u\in L^{m+1}(\Omega).$ This ends the proof of the lemma.
\medskip
\end{proof*}

\begin{lem}
\label{lemdualpro}
Let $u_1,u_2\in H^1_0(\Omega),$ $p\in[1,\infty)$ and $v_1,v_2\in L^{p^\p}(\Omega)$ be such that $\Delta u_j+v_j\in L^2(\Omega),$ for any $j\in\{1,2\}.$ We then have,
\begin{align*}
	&	\; \big((\Delta u_1+v_1)-(\Delta u_2+v_2),w_1-w_2\big)_{L^2(\Omega)}		\\
   =	&	\; -\big(\nabla(u_1-u_2),\nabla(w_1-w_2)\big)_{L^2(\Omega)}+\langle v_1-v_2,w_1-w_2\rangle_{L^{p^\p}(\Omega),L^p(\Omega)},
\end{align*}
for any $w_1,w_2\in H^1_0(\Omega)\cap L^p(\Omega).$
\end{lem}

\begin{proof*}
Let $X=H^1_0(\Omega)\cap L^p(\Omega).$ We recall that since $H^1_0(\Omega)\cap L^p(\Omega)$ is dense in both $H^1_0(\Omega)$ and $L^p(\Omega),$ we have by \eqref{dual} that $X^\star=H^{-1}(\Omega)+L^{p^\p}(\Omega).$ We also recall that we identify $L^2(\Omega)$ with its own dual, so that, by \eqref{dualg}, the $L^2$-scalar product is also the $L^2-L^2$ duality product. Finally, by \eqref{dualtranspose}, since the embeddings of $X$ in $L^2(\Omega),$ $L^p(\Omega)$ and $H^1_0(\Omega)$ are all continuous and dense, it follows that for any $j\in\{1,2\},$ $\Delta u_j,v_j\in X^\star$ and
\begin{align*}
	&	\; \big((\Delta u_1+v_1)-(\Delta u_2+v_2),w_1-w_2\big)_{L^2(\Omega)}
		=\langle(\Delta u_1+v_1)-(\Delta u_2+v_2),w_1-w_2\rangle_{X^\star,X}						\\
   =	&	\; \langle\Delta(u_1-u_2),w_1-w_2\rangle_{X^\star,X}+\langle v_1-v_2,w_1-w_2\rangle_{X^\star,X},	\\
   =	&	\; \langle\Delta(u_1-u_2),w_1-w_2\rangle_{H^{-1}(\Omega),H^1_0(\Omega)}
			+\langle v_1-v_2,w_1-w_2\rangle_{L^{p^\p}(\Omega),L^p(\Omega)},
\end{align*}
from which we deduce the result.
\medskip
\end{proof*}

\begin{cor}
\label{corAmon}
Assume~\eqref{O}--\eqref{a}. The operator $(A,D(A))$ is monotone on $L^2(\Omega).$
\end{cor}

\begin{proof*}
Let $u,v\in D(A).$ Let $X=H^1_0(\Omega)\cap L^{m+1}(\Omega).$ It follows from Lemma~\ref{lemdualpro}, \eqref{dualg} and \cite[Corollary~5.8]{MR4340780} that,
\begin{gather*}
(Au-Av,u-v)_{L^2(\Omega)}=\langle-\vi a(g(u)-g(v)),u-v\rangle_{L^\frac{m+1}m(\Omega),L^{m+1}(\Omega)}\ge0.
\end{gather*}
Hence the result.
\medskip
\end{proof*}

\begin{vproof}{of Lemma~\ref{lemA}.}
The density is obvious. By Corollary~\ref{corAmon} and Brezis~\cite[Proposition~2.2]{MR0348562}, we only have to show that $R(I+A)=L^2(\Omega).$ Let $F\in L^2(\Omega).$ Let $u\in H^1_0(\Omega)\cap L^{m+1}(\Omega)$ and $u_{\eps_n}\in D(A_{\eps_n}^m)$ $(n\in\N)$ be given by Lemma~\ref{lemAmax}. It follows from \eqref{lemAmax2} and \eqref{lemAmax6} that $\big(g_{\eps_n}^m(u_{\eps_n})\big)_{n\in\N}$ is bounded in $L^\frac2m(\Omega)$ and that $g_{\eps_n}^m(u_{\eps_n})\xrightarrow[n\to\infty]{\text{a.e.\,in }\Omega}g(u).$ Thus by Strauss~\cite{MR306715},
\begin{gather}
\label{demlemA1}
g_{\eps_n}^m(u_{\eps_n})\xrightarrow[n\to\infty]{\Dr^\p(\Omega)}g(u).
\end{gather}
Passing to the limit as $n\tends\infty$ in \eqref{lemAmax1}, it follows from \eqref{lemAmax4}, \eqref{lemAmax5} and \eqref{demlemA1} that $u$ satisfies
\begin{gather}
\label{demlemA2}
-\vi\Delta u-\vi Vu-\vi ag(u)+u=F, \text{ in } \Dr^\p(\Omega).
\end{gather}
But $u\in H^1_0(\Omega)\cap L^{m+1}(\Omega)$ and $Vu,F\in L^2(\Omega)$ so that, by \eqref{demlemA2},
\begin{gather*}
u\in D(A) \; \text{ and } \; u+Au=F, \text{ in } L^2(\Omega).
\end{gather*}
This concludes the proof.
\medskip
\end{vproof}

\begin{vproof}{of Proposition~\ref{propsolL2}.}
Set $Y=H^2_0(\Omega)\cap L^\frac2{2-m}(\Omega).$ By~\eqref{dual}, $Y^\star=H^{-2}(\Omega)+L^\frac2m(\Omega).$ By \eqref{cv}, \eqref{g} and \eqref{lemVH-1}, we have for any $T>0,$
\begin{align}
\label{dempropsolL21}
&	\Delta u_n\xrightarrow[n\to\infty]{C([0,T];H^{-2}(\Omega))}\Delta u,		\\
\label{dempropsolL22}
&	Vu_n\xrightarrow[n\to\infty]{C([0,T];H^{-1}(\Omega))}Vu,				\\
\label{dempropsolL23}
&	g(u_n)\xrightarrow[n\to\infty]{C([0,T];L^\frac2m(\Omega))}g(u),
\end{align}
Then it follows from the equation satisfied by each $u_n,$ \eqref{cv} and \eqref{dempropsolL21}--\eqref{dempropsolL23} that for any $T>0,$ $(u_n)_{n\in\N}$ is a Cauchy sequence in $W^{1,1}\big((0,T);Y^\star\big),$ so that \eqref{propsolL21}--\eqref{propsolL22} hold true. We use \eqref{cv}, \eqref{propsolL22} and \eqref{dempropsolL21}--\eqref{dempropsolL23} to pass in the limit in the equation satisfied by each $u_n.$ With help of \eqref{injDp}, it follows that $u$ satisfies \eqref{nls} in $L^1\big([0,\infty);Y^\star\big)\inj\Dr^\p\big((0,\infty)\times\Omega\big).$
\medskip
\end{vproof}

\begin{vproof}{of Proposition~\ref{propdep}.}
By \cite[Proposition~2.5]{MR4340780}, we only have to show \eqref{estthmweak} for the weak solutions. The $H^1_0$-solutions satisfying \eqref{propdephyp} with $p=m+1,$ and estimate~\eqref{estthmweak} being stable by passing to the limit in $L^1\big((0,T);L^2(\Omega)\big)\times C\big([0,T];L^2(\Omega)\big),$ the result is then obtained with help of \eqref{cv}.
\medskip
\end{vproof}

\begin{vproof}{of Theorem~\ref{thmstrongaH2}.}
Let $f$ and $u_0$ be as in the theorem. By Lemma~\ref{lemA} and Barbu~\cite[Theorem~4.5, p.141]{MR2582280} (see also Vrabie~\cite[Theorem~1.7.1, p.23]{MR1375237}), there exists a unique solution $u\in W^{1,\infty}_\loc\big([0,\infty);L^2(\Omega)\big)$ to \eqref{nls}--\eqref{u0} satisfying for almost every $t>0,$ $u(t)\in D(A)$ and \eqref{strongaH23}, from which \eqref{strongaH21} follows. Since $u\in W^{1,\infty}_\loc\big([0,\infty);L^2(\Omega)\big),$ it follows from Lemma~A.5 in~Bégout and D\'iaz~\cite{MR4053613} that $M:t\longmapsto\frac12\|u(t)\|_{L^2(\Omega)}^2$ belongs to $W^{1,\infty}_\loc\big([0,\infty);\R\big)$ and $M^\p(t)=\big(u_t(t),u(t)\big)_{L^2(\Omega)},$ for almost every $t>0.$ Taking the $L^2$-scalar product of \eqref{nls} with $\vi u,$ we get Property~\ref{thmstrongaH23}, with help of Lemma~\ref{lemdualpro}. We apply the Cauchy-Schwarz inequality to \eqref{L2} and we inegrate in time to obtain \eqref{strongaH22}. Now, we take again the $L^2$-scalar product of \eqref{nls} with $-u.$ We use Lemma~\ref{lemdualpro} and the fact that $a\in D(m).$ We sum the result with $\left(\frac{2\sqrt m}{1-m}+1\right)\times\eqref{L2}.$ Finally, we use again the Cauchy-Schwarz inequality to infer that,
\begin{gather*}
\|\nabla u\|_{L^2(\Omega)}^2+\Im(a)\|u\|_{L^{m+1}(\Omega)}^{m+1}
\le C(m)\left(\|u_t\|_{L^2(\Omega)}+\|Vu\|_{L^2(\Omega)}+\|f\|_{L^2(\Omega)}\right)\|u\|_{L^2(\Omega)},
\end{gather*}
almost everywhere on $(0,\infty).$ It follows from \eqref{lemV1}--\eqref{lemV2} that,
\begin{gather}
\label{demlemthmsaH21b}
\begin{split}
	&	\; \|\nabla u\|_{L^2(\Omega)}^2+\Im(a)\|u\|_{L^{m+1}(\Omega)}^{m+1}					\\
  \le	&	\; C(N,m,\beta)\left(\|u_t\|_{L^2(\Omega)}+\left(\|V_1\|_{L^\infty(\Omega)}+\|V_2\|_{L^{p_V}(\Omega)}^{2-\gamma}\right)
  			\|u\|_{L^2(\Omega)}+\|f\|_{L^2(\Omega)}\right)\|u\|_{L^2(\Omega)},
\end{split}
\end{gather}
from which \eqref{strongaH24} follows. Then Property~\ref{thmstrongaH22} holds true, from which we deduce Property~\ref{thmstrongaH24}. Moreover, since $u\in C\big([0,\infty);L^2(\Omega)\big),$ Property~\ref{thmstrongaH21} comes from \eqref{strongaH24} and \eqref{rmkweakcon}. Finally, it follows from Property~\ref{thmstrongaH21} that $u$ is an $H^1_0$-solution and that $u$ satisfies~\eqref{nls} in $L^\infty_\loc\big([0,\infty);H^{-1}(\Omega)+L^\frac{m+1}m(\Omega)\big).$ The theorem is proved.
\medskip
\end{vproof}

\begin{vproof}{of Theorem~\ref{thmweak}.}
By Theorem~\ref{thmstrongaH2}, Proposition~\ref{propdep} and \ref{rmkdefsol1} of Remark~\ref{rmkdefsol}, the proof follows easily by density of $\Dr(\Omega)\times W^{1,1}_\loc\big([0,\infty);L^2(\Omega)\big)$ in $L^2(\Omega)\times L^1_\loc\big([0,\infty);L^2(\Omega)\big)$ (see the proof of Theorem~2.6 in \cite{MR4340780} for more details).
\medskip
\end{vproof}

\noindent
We split the proof of Theorems~\ref{thmstrongaH1} and \ref{thmstrongH1} into several lemmas.

\begin{lem}
\label{lemthmsaH1}
Let Assumption~$\ref{ass1}$ be fulfilled with additionally $V\in W^{1,\infty}(\Omega;\R)+W^{1,p_V}(\Omega;\R).$ Let $f$ satisfy \eqref{fH1} and $u_0\in H^1_0(\Omega).$ Let $(f_\eps)_{\eps>0}\subset\Dr\big([0,\infty);H^1_0(\Omega)\big)$ and $(\vphi_\eps)_{\eps>0}\subset\Dr(\Omega)$ be such that,
\begin{gather}
\label{lemthmsaH11}
\begin{cases}
f_\eps\xrightarrow[\eps\searrow0]{L^1((0,T);H^1_0(\Omega))\cap L^\frac{m+1}m((0,T);X^\star)}f,		\\
\vphi_\eps\xrightarrow[\eps\searrow0]{H^1_0(\Omega)}u_0.
\end{cases}
\end{gather}
for any $T>0,$ where $X=H^1_0(\Omega)\cap L^{m+1}(\Omega).$ Then for any $\eps>0,$ there exists a unique solution
\begin{gather}
\label{lemthmsaH12}
u_\eps\in C_\w\big([0,\infty);H^1_0(\Omega)\big)\cap W^{1,\infty}_\loc\big([0,\infty);L^2(\Omega)\big),
\end{gather}
to,
\begin{gather}
\label{nlse}
\vi\frac{\partial u_\eps}{\partial t}+\Delta u_\eps+V(x)u_\eps+ag_\eps^m(u_\eps)=f_\eps(t,x),	\text{ in } L^2(\Omega),
\end{gather}
such that $u_\eps(0)=\vphi_\eps.$ Furthermore, the following holds for any $\eps>0.$
\begin{gather}
\label{lemthmsaH13}
\|u_\eps(t)\|_{H^1_0(\Omega)}
\le\left(\|\vphi_\eps\|_{H^1_0(\Omega)}+\vint_0^t\|f_\eps(\sigma)\|_{H^1_0(\Omega)}\d\sigma\right)e^{C\|\nabla V\|_{L^\infty+L^{p_V}}t},
\end{gather}
for any $t\ge0,$ where $C=C(N,\beta),$ and if $\nabla V=0$ then,
\begin{gather}
\label{lemthmsaH14}
\|\nabla u_\eps(t)\|_{L^2(\Omega)}\le\|\nabla\vphi_\eps\|_{L^2(\Omega)}+\vint_0^t\|\nabla f_\eps(\sigma)\|_{L^2(\Omega)}\d\sigma,
\end{gather}
for any $t\ge0.$ Finally,
\begin{gather}
\label{lemthmsaH15}
\begin{cases}
(u_\eps)_{\eps>0} \text{ is bounded in }
	L^\infty_\loc\big([0,\infty);H^1_0(\Omega)\big)\cap W^{1,\frac{m+1}m}_\loc\big([0,\infty);X^\star+L^\frac2m(\Omega)\big),	\medskip \\
\big(g_\eps^m(u_\eps)\big)_{\eps>0} \text{ is bounded in } L^\infty_\loc\big([0,\infty);L^\frac2m(\Omega)\big),
\end{cases}
\end{gather}
and
\begin{gather}
\label{lemthmsaH16}
\sup_{\eps>0}\vint_0^T\!\!\!\vint_\Omega\frac{|u_\eps(t,x)|^2}{(|u_\eps(t,x)|^2+\eps)^\frac{1-m}2}\d x\d t\le C(T),
\end{gather}
for any $T>0.$
\end{lem}

\begin{proof*}
Let the assumptions of the Lemma be fulfilled. By Lemma~\ref{lemAe} and Barbu~\cite[Theorem~4.5, p.141]{MR2582280} (see also Vrabie~\cite[Theorem~1.7.1, p.23]{MR1375237}), there exists a unique solution $u_\eps\in W^{1,\infty}_\loc\big([0,\infty);L^2(\Omega)\big)$ to \eqref{nlse} such that $u_\eps(0)=\vphi_\eps.$ Moreover, $u_\eps(t)\in D(A_\eps^m),$ for almost every $t>0.$ Now, we take the $L^2$-scalar product of \eqref{nlse} with $-u_\eps$ and we get with help of Cauchy-Schwarz's inequality that,
\begin{gather*}
\|\nabla u_\eps\|_{L^2(\Omega)}^2
\le\left(\|u_\eps^\p\|_{L^2(\Omega)}+\|Vu_\eps\|_{L^2(\Omega)}+|a|\eps^{-\frac{1-m}2}\|u_\eps\|_{L^2(\Omega)}+\|f_\eps\|_{L^2(\Omega)}\right)\|u_\eps\|_{L^2(\Omega)},
\end{gather*}
almost everywhere on $(0,\infty).$ Applying \eqref{lemV1} and \eqref{lemV2} to the above with $\rho=2\|u_\eps\|_{L^2(\Omega)},$ we get that $u_\eps\in L^\infty_\loc\big([0,\infty);H^1_0(\Omega)\big).$ With help of \eqref{rmkweakcon}, \eqref{lemthmsaH12} follows. By \eqref{lemVL2} and \eqref{nlse}, it follows that $\Delta u_\eps\in L^\infty_\loc\big([0,\infty);H^1_0(\Omega)\big).$ So, we are allowed to apply \cite[Lemma~A.5]{MR4053613}. Taking the $L^2$-scalar product of \eqref{nlse} with $-\vi\Delta u_\eps,$ it then follows from~(6.8) in \cite{MR4053613} and a density argument that for almost every $\sigma>0,$
\begin{gather}
\label{demlemthmsaH11}
\frac12\frac{\d}{\d t}\|\nabla u_\eps(\sigma)\|_{L^2(\Omega)}^2\le\big(\nabla f_\eps(\sigma)-u_\eps(\sigma)\nabla V,\vi\nabla u_\eps(\sigma)\big)_{L^2(\Omega)}.
\end{gather}
Let $t>0.$ If $\nabla V=0$ then it follows from \eqref{demlemthmsaH11} and Cauchy-Schwarz's inequality that,
\begin{gather*}
\frac12\frac{\d}{\d t}\|\nabla u_\eps(\sigma)\|_{L^2(\Omega)}^2\le\|\nabla f_\eps(\sigma)\|_{L^2(\Omega)}\|\nabla u_\eps(\sigma)\|_{L^2(\Omega)}.
\end{gather*}
Integrating over $(0,t),$ we obtain \eqref{lemthmsaH14}. Now, we turn out to the general case. Taking the $L^2$-scalar product of \eqref{nlse} with $\vi u_\eps,$ we get with help of Bégout and D\'iaz~\cite[Lemma~A.5]{MR4053613} and Cauchy-Schwarz's inequality that,
\begin{gather}
\label{demlemthmsaH12}
\frac12\frac{\d}{\d t}\|u_\eps(\sigma)\|_{L^2(\Omega)}^2+\Im(a)\vint_\Omega\frac{|u_\eps(\sigma,x)|^2}{(|u_\eps(\sigma,x)|^2+\eps)^\frac{1-m}2}\d x
\le\|f_\eps(\sigma)\|_{L^2(\Omega)}\|u_\eps(\sigma)\|_{L^2(\Omega)},
\end{gather}
for almost every $\sigma>0.$ Now, let us still apply Cauchy-Schwarz's inequality in \eqref{demlemthmsaH11}. Using \eqref{lemVL2} and summing the result with \eqref{demlemthmsaH12}, we get for almost every $\sigma>0,$
\begin{gather*}
\frac12\frac{\d}{\d t}\|u_\eps(\sigma)\|_{H^1_0(\Omega)}^2
\le\|f_\eps(\sigma)\|_{H^1_0(\Omega)}\|u_\eps(\sigma)\|_{H^1_0(\Omega)}+C\|\nabla V\|_{L^\infty(\Omega)+L^{p_V}(\Omega)}\|u_\eps(\sigma)\|_{H^1_0(\Omega)}^2,
\end{gather*}
where $C$ is given by \eqref{lemVL2}. Integrating over $(0,t),$ we obtain
\begin{gather*}
\|u_\eps(t)\|_{H^1_0}\le\underbrace{\|u_0\|_{H^1_0}+\vint_0^t\|f(\sigma)\|_{H^1_0}\d\sigma}_{\|\atop\vphi(t)}
+\vint_0^t\underbrace{C\|\nabla V\|_{L^\infty+L^{p_V}}}_{\|\atop\alpha}\|u_\eps(\sigma)\|_{H^1_0}\d\sigma,
\end{gather*}
and by Gronwall's Lemma (see, for instance, Barbu~\cite[Lemma~1.1, p.22]{MR3585801}),
\begin{align*}
\|u_\eps(t)\|_{H^1_0}		&	\le\vphi(t)+\vint_0^t\alpha\vphi(\sigma)\exp\left(\vint_\sigma^t\alpha\d s\right)\d\sigma
										=\vphi(t)+\vint_0^t\alpha\vphi(\sigma)e^{\alpha(t-\sigma)}\d\sigma	\\
					&	\le\vphi(t)+\vphi(t)\vint_0^t\alpha e^{\alpha(t-\sigma)}\d\sigma=\vphi(t)e^{\alpha t},
\end{align*}
which is \eqref{lemthmsaH13}. And from \eqref{lemthmsaH13} we get that $(u_\eps)_{\eps>0}$ is bounded in $L^\infty_\loc\big([0,\infty);H^1_0(\Omega)\big),$ from which we deduce the boundness of $\big(g_\eps^m(u_\eps)\big)_{\eps>0}$ in $L^\infty_\loc\big([0,\infty);L^\frac2m(\Omega)\big),$ since for any $\eps>0,$ $|g_\eps^m(u_\eps)|\le|u_\eps|^m,$ almost everywhere in $(0,\infty)\times\Omega.$ It then follows from \eqref{lemthmsaH11}, and the equation \eqref{nlse} that $(u_\eps)_{\eps>0}$ is bounded in $W^{1,\frac2m}_\loc\big([0,\infty);X^\star+L^\frac2m(\Omega)\big),$ from which \eqref{lemthmsaH15} follows. Let $T>0.$ Integrating \eqref{demlemthmsaH12} over $(0,T),$ applying Hölder's inequality, and using \eqref{lemthmsaH11} and \eqref{lemthmsaH15}, we get \eqref{lemthmsaH16}. The lemma is proved.
\medskip
\end{proof*}

\begin{lem}
\label{lemthmstrongaH23}
Let Assumption~$\ref{ass1}$ be fulfilled with additionally $V\in W^{1,\infty}(\Omega;\R)+W^{1,p_V}(\Omega;\R).$ We use the notations of Lemma~$\ref{lemthmsaH1}.$ Under the hypotheses of Lemma~$\ref{lemthmsaH1},$ there exist
\begin{gather}
\label{lemthmstrongaH23a}
u\in C_\w\big([0,\infty);H^1_0(\Omega)\big)\cap W^{1,\frac{m+1}m}_\loc\big([0,\infty);X^\star+L^\frac2m(\Omega)\big),	\\
\label{lemthmstrongaH23b}
u\in  L^{m+1}_\loc\big([0,\infty);L^{m+1}(\Omega)\big),
\end{gather}
and a positive sequence $\eps_n\searrow0,$ as $n\tends\infty,$ such that
\begin{align}
\label{lemthmstrongaH23c}
&	u_{\eps_n}(t)\underset{n\to\infty}{-\!\!\!-\!\!\!-\!\!\!-\!\!\!\weak}u(t) \text{ in } H^1_0(\Omega)_\w, \; \forall t\ge0,		\\
\label{lemthmstrongaH23d}
&	u_{\eps_n}\xrightarrow[n\to\infty]{\text{a.e.\,in }(0,\infty)\times\Omega}u,									\\
\label{lemthmstrongaH23e}
&	g_{\eps_n}^m\big(u_{\eps_n}\big)\underset{n\to\infty}{-\!\!\!-\!\!\!-\!\!\!-\!\!\!\weak}g(u) \text{ in } L^\frac2m\big((0,T);L^\frac2m(\Omega)\big)_\w,
\end{align}
for any $T>0.$
\end{lem}

\begin{proof*}
We first note that,
\begin{gather}
\label{demlemthmsaH23a}
W^{1,\frac{m+1}m}_\loc\big([0,\infty);X^\star+L^\frac2m(\Omega)\big)\inj C^{0,\frac1{m+1}}_\loc\big([0,\infty);X^\star+L^\frac2m(\Omega)\big).
\end{gather}
By \eqref{lemthmsaH15}, \eqref{demlemthmsaH23a}, Cazenave~\cite{MR2002047} (Proposition~1.1.2(i), p.2, and Remark~1.3.13(ii), p.12) and the diagonal procedure, we obtain the existence of a
\begin{gather}
\label{demlemthmsaH23b}
u\in C_\w\big([0,\infty);H^1_0(\Omega)\big)\cap W^{1,\frac{m+1}m}_\loc\big([0,\infty);X^\star+L^\frac2m(\Omega)\big)
\end{gather}
satisfying \eqref{lemthmstrongaH23c}. Let $T>0$ and $\Omega^\p\subset\Omega$ be any bounded open subset of $\R^N$ having a $C^1$-boundary. By Rellich-Kondrachov's compactness Theorem, we have that,
\begin{gather}
\label{demlemthmsaH23c}
H^1(\Omega^\p)\underset{\text{compact}}{\inj} L^2(\Omega^\p)\inj H^{-1}(\Omega^\p)+L^\frac{m+1}m(\Omega^\p)+L^\frac2m(\Omega^\p),
\end{gather}
and by \eqref{lemthmsaH15},
\begin{gather}
\left\{
\begin{split}
\label{demlemthmsaH23d}
(u_\eps)_{\eps>0} \text{ is bounded in } L^\infty_\loc\big([0,\infty);H^1(\Omega^\p)\big) \text{ and in}	\\
W^{1,\frac{m+1}m}_\loc\big([0,\infty);H^{-1}(\Omega^\p)+L^\frac{m+1}m(\Omega^\p)+L^\frac2m(\Omega^\p)\big).
\end{split}
\right.
\end{gather}
It follows from \eqref{lemthmstrongaH23c}, \eqref{demlemthmsaH23c}--\eqref{demlemthmsaH23d} and a compactness result due to Simon~\cite{MR916688} (Corollary~4, p.85) that,
\begin{gather*}
u\in C\big([0,T];L^2(\Omega^\p)\big) \text{ and } \lim_{n\to\infty}\|u_{\eps_n}-u\|_{C([0,T];L^2(\Omega^\p))}=0.
\end{gather*}
Since $T$ and $\Omega^\p$ are arbitrary, we deduce that $u_{\eps_n}\xrightarrow[n\to\infty]{L^2_\loc((0,\infty)\times\Omega)}u.$ Up to a subsequence, that we still denote by $(u_{\eps_n})_{n\in\N},$ and with help of the diagonal procedure, we obtain \eqref{lemthmstrongaH23d}. It follows from \eqref{lemthmstrongaH23d} that,
\begin{gather}
\label{demlemthmsaH23e}
g_{\eps_n}^m\big(u_{\eps_n}\big)\xrightarrow[n\to\infty]{\text{a.e.\,in }(0,\infty)\times\Omega}g(u).
\end{gather}
Since for any $T>0,$ $L^\frac2m\big((0,T);L^\frac2m(\Omega)\big)\cong L^\frac2m\big((0,T)\times\Omega\big),$ we obtain \eqref{lemthmstrongaH23e} from \eqref{lemthmsaH15}, \eqref{demlemthmsaH23e}, and Cazenave~\cite[Proposition~1.2.1, p.3]{MR2002047}. Finally, \eqref{lemthmstrongaH23b} is a consequence of \eqref{lemthmsaH16}, \eqref{lemthmstrongaH23d} and Fatou's Lemma.
\medskip
\end{proof*}

\begin{lem}
\label{lemthmstrongaH24}
Let Assumption~$\ref{ass1}$ be fulfilled with additionally $V\in W^{1,\infty}(\Omega;\R)+W^{1,p_V}(\Omega;\R).$ We use the notations of Lemma~$\ref{lemthmsaH1}.$ Under the hypotheses of Lemma~$\ref{lemthmsaH1},$ the function $u$ given by Lemma~$\ref{lemthmstrongaH23}$ is the unique $H^1_0$-solution to \eqref{nls}--\eqref{u0}. In addition, $u$ satisfies \eqref{thmstrongaH12} and \eqref{thmstrongaH13} with $s=0,$ according to the different cases satisfied by $V.$
\end{lem}

\begin{proof*}
Let $u$ be given by Lemma~\ref{lemthmstrongaH23}. Uniqueness comes from Proposition~\ref{propdep}. Set $Y=L^\frac2{2-m}(\Omega).$ Let $\vphi\in X\cap Y$ and $\psi\in C^1_\co\big((0,\infty);\R\big).$ Let $T>0$ be such that $\supp\psi\in(0,T).$ Let $(\eps_n)_{n\in\N}$ be given by Lemma~\ref{lemthmstrongaH23}. It follows from \eqref{nlse} and \eqref{lemdualV} that for any $n\in\N,$
\begin{gather*}
\vint_0^\infty\left\langle\vi\frac{\partial u_{\eps_n}}{\partial t}+\Delta u_{\eps_n}+Vu_{\eps_n}
	+ag_{\eps_n}^m\big(u_{\eps_n}\big),\vphi\right\rangle_{X^\star+Y^\star,X\cap Y}\psi(t)\,\d t=\vint_0^\infty\big\langle f_{\eps_n}(t),\vphi\big\rangle_{X^\star+Y^\star,X\cap Y}\psi(t)\,\d t,
\end{gather*}
and so,
\begin{align*}
	&	\; \vint_0^T\Big(\left\langle-\vi u_{\eps_n},\vphi\right\rangle_{L^2(\Omega),L^2(\Omega)}\psi^\p(t)
			-\left\langle\nabla u_{\eps_n},\nabla\vphi\right\rangle_{L^2(\Omega),L^2(\Omega)}\psi(t)
			+\left\langle u_{\eps_n},V\vphi\right\rangle_{L^2(\Omega),L^2(\Omega)}\psi(t)									\\
   +	&	\; \left\langle ag_{\eps_n}^m\big(u_{\eps_n}\big),\vphi\right\rangle_{Y^\star,Y}\psi(t)\Big)\d t
			=\vint_0^T\big\langle f_{\eps_n}(t),\vphi\big\rangle_{X^\star,X}\psi(t)\,\d t.
\end{align*}
By \eqref{lemthmsaH11}, \eqref{lemthmstrongaH23c}, \eqref{lemthmstrongaH23e} and the dominated convergence Theorem, we can pass to the limit in the above equality to obtain $u(0)=u_0,$ and
\begin{gather*}
\vint_0^\infty\left\langle\vi\frac{\partial u}{\partial t}+\Delta u+Vu+ag(u),\vphi\right\rangle_{X^\star+Y^\star,X\cap Y}\psi(t)\,\d t
=\vint_0^\infty\big\langle f(t),\vphi\big\rangle_{X^\star+Y^\star,X\cap Y}\psi(t)\,\d t.
\end{gather*}
It follows that $u$ satisfies \eqref{nls} in $L^1_\loc\big((0,\infty);X^\star+Y^\star\big),$ hence in $\Dr^\p\big((0,\infty)\times\Omega\big).$ From \eqref{nls}, \eqref{lemthmstrongaH23a}, \eqref{lemthmstrongaH23b}, and \eqref{lemVL2}, we deduce that $u\in W^{1,\frac{m+1}m}_\loc\big([0,\infty);X^\star\big).$ So, $u$ is the unique $H^1_0$-solution. Finally, \eqref{thmstrongaH12} and \eqref{thmstrongaH13} for $s=0$ come from \eqref{lemthmsaH11}, \eqref{lemthmsaH13}, \eqref{lemthmsaH14}, \eqref{lemthmstrongaH23c} and the lower semicontinuity of the norm. This ends the proof of the lemma.
\medskip
\end{proof*}

\begin{vproof}{of Theorem~\ref{thmstrongH1}.}
Let $u_0\in H^1_0(\Omega)$ and let $f$ satisfy \eqref{fH1}. Let $u$ be given by Lemma~\ref{lemthmstrongaH23}. By Lemma~\ref{lemthmstrongaH24}, it remains to show that $u$ satisfies \eqref{L2}. Let $X=H^1_0(\Omega)\cap L^{m+1}(\Omega).$ Taking the $X^\star-X$ duality product of \eqref{nls} with $\vi u,$ and applying~\cite[Lemma~A.5]{MR4053613} and \eqref{dualg}, we obtain \eqref{L2}. This ends the proof of Theorem~\ref{thmstrongH1}.
\medskip
\end{vproof}

\begin{vproof}{of Theorems~\ref{thmstrongaH1}.}
Let $u_0\in H^1_0(\Omega)$ and $f\in L^1_\loc\big([0,\infty);H^1_0(\Omega)\big).$ Let $(\vphi_\eps)_{\eps>0}\subset\Dr(\Omega)$ and $(f_n)_{n\in\N}\subset\Dr\big([0,\infty);H^1_0(\Omega)\big)$ be such that $\vphi_n\xrightarrow[n\to\infty]{H^1_0(\Omega)}u_0$ and $f_n\xrightarrow[n\to\infty]{L^1((0,T);H^1_0)}f,$ for any $T>0.$ For each $n\in\N,$ let $u_n$ be the unique $H^1_0$-solution to \eqref{nls} such that $u_n(0)=\vphi_n,$ given by Lemma~\ref{lemthmstrongaH24}. By Proposition~\ref{propdep}, $(u_n)_{n\in\N}$ is a Cauchy sequence in $C\big([0,T];L^2(\Omega)\big),$ for any $T>0.$ As a consequence, there exists $u\in C\big([0,\infty);L^2(\Omega)\big)$ such that for any $T>0,$
\begin{gather}
\label{demthmstrongaH12}
u_n\xrightarrow[n\to\infty]{C([0,T];L^2(\Omega))}u.
\end{gather}
By definition, $u$ is a weak solution and satisfies \eqref{nls} in $\Dr^\p\big((0,\infty)\times\Omega\big)$ (Proposition~\ref{propsolL2}). In particular, $u$ fulfills \eqref{Lm}. Still by Lemma~\ref{lemthmstrongaH24}, each $u_n$ satisfies \eqref{thmstrongaH12} and \eqref{thmstrongaH13} with $s=0$ so that,
\begin{gather}
\label{demthmstrongaH14}
(u_n)_{n\in\N} \text{ is bounded in } L^\infty_\loc\big([0,\infty);H^1_0(\Omega)\big).
\end{gather}
We deduce from \eqref{rmkweakcon}, \eqref{demthmstrongaH12} and \eqref{demthmstrongaH14} that $u\in C_\w\big([0,\infty);H^1_0(\Omega)\big)$ and,
\begin{gather}
\label{demthmstrongaH15}
\forall t\ge0, \; u_n(t)\underset{n\to\infty}{-\!\!\!-\!\!\!-\!\!\!-\!\!\!\weak}u(t) \text{ in } H^1_0(\Omega)_\w,
\end{gather}
Then $u$ satisfies the first line of \eqref{thmstrongaH11}. By \eqref{nls}, \eqref{Lm}, the first line of \eqref{thmstrongaH11}, \eqref{g} and \eqref{lemVL2}, $u$ satisfies the second line of \eqref{thmstrongaH11}, and \eqref{nls} in $L^1_\loc\big([0,\infty);H^{-1}(\Omega)+L^\frac{m+1}m(\Omega)\big).$ Passing to the limit, as $n\tends\infty,$ in \eqref{thmstrongaH12}--\eqref{thmstrongaH13} satisfied by each $u_n,$ and using \eqref{demthmstrongaH15} and the lower semicontinuity of the norm, we obtain \eqref{thmstrongaH12}--\eqref{thmstrongaH13} for $u$ with $s=0.$ The general case follows by standard arguments of time translation and uniqueness of the weak solutions. See, for instance, the end of the proof of \cite[Theorem~2.7]{MR4340780}.
\medskip
\end{vproof}

\section{On the $\bs{H^2}$-solutions}
\label{H2sol}

\begin{defi}
\label{defsolH2}
Assume~\eqref{O}, \eqref{m}, \eqref{V} and \eqref{pV}. Let $a\in\C,$ $f\in L^1_\loc\big([0,\infty);L^2(\Omega)\big)$ and $u_0\in L^2(\Omega).$ We shall say that $u$ is an $H^2$-\textit{solution of} \eqref{nls}--\eqref{u0} if $u$ is an $H^1_0$-\textit{solution of} \eqref{nls}--\eqref{u0}, if $u\in W^{1,\frac{m+1}m}_\loc\big([0,\infty);L^2(\Omega)+L^\frac{m+1}m(\Omega)\big)$ and if for almost every $t>0,$ $\Delta u(t)\in L^2(\Omega).$
\end{defi}

\begin{thm}[\textbf{Existence and uniqueness of $\bs{H^2}$-solutions}]
\label{thmstrongH2}
Let Assumption~$\ref{ass1}$ be fulfilled and $f\in W^{1,1}_\loc\big([0,\infty);L^2(\Omega)\big).$ If $|\Omega|<\infty$ then for any $u_0\in H^1_0(\Omega)$ for which $\Delta u_0\in L^2(\Omega),$ there exists a unique $H^2$-solution $u$ to \eqref{nls}--\eqref{u0}. Furthermore, $u$ satisfies~\eqref{nls} in $L^\infty_\loc\big([0,\infty);L^2(\Omega)\big)$ as well as the following properties.
\begin{enumerate}
\item
\label{thmstrongH21}
$u\in C\big([0,\infty);H^1_0(\Omega)\big)\cap W^{1,\infty}_\loc\big([0,\infty);L^2(\Omega)\big).$
\item
\label{thmstrongH22}
$\Delta u\in C_\w\big([0,\infty);L^2(\Omega)\big)$ and for any $t\ge s\ge0,$
\begin{gather}
\label{thmstrongH22a}
\|\nabla u(t)-\nabla u(s)\|_{L^2(\Omega)}\le2\|u_t\|_{L^\infty((s,t);L^2(\Omega))}^\frac12\|\Delta u\|_{L^\infty((s,t);L^2(\Omega))}^\frac12|t-s|^\frac12.
\end{gather}
\item
\label{thmstrongH23}
The map $t\longmapsto\|u(t)\|_{L^2(\Omega)}^2$ belongs to $C^1\big([0,\infty);\R\big)$  and~\eqref{L2} holds for any $t\ge0.$
\item
\label{thmstrongH24}
If $f\in W^{1,1}\big((0,\infty);L^2(\Omega)\big)$ then we have,
\begin{align*}
	&	u\in C_\b\big([0,\infty);H^1_0(\Omega)\big)\cap W^{1,\infty}\big((0,\infty);L^2(\Omega)\big),	\\
	&	\Delta u\in L^\infty\big((0,\infty);L^2(\Omega)\big).
\end{align*}
\end{enumerate}
\end{thm}

\begin{proof*}
Let $f$ and $u$ be as in the theorem. Since $|\Omega|<\infty,$ we have that $u_0\in L^{m+1}(\Omega)$ and $g(u_0)\in L^2(\Omega).$ It follows that Theorem~\ref{thmstrongaH2} applies. It follows easily from \eqref{nls} that $u,$ which is given by Theorem~\ref{thmstrongaH2}, satisfies,
\begin{gather}
\label{demthmsH21}
\Delta u\in L^\infty_\loc\big([0,\infty);L^2(\Omega)\big),
\end{gather}
and \eqref{nls} makes sense in $L^\infty_\loc\big([0,\infty);L^2(\Omega)\big).$ As a consequence, $u$ is an $H^2$-solution. But any $H^2$-solution is an $H^1_0$-solution, for which we have uniqueness, so that $u$ is the unique solution. Since $\Delta u\in C\big([0,\infty);H^{-2}(\Omega)\big)$ and $u\in C\big([0,\infty);L^2(\Omega)\big),$ Properties~\ref{thmstrongH21}--\ref{thmstrongH23} are then obtained from \eqref{rmkweakcon}, \eqref{demthmsH21}, \eqref{nablau} and Properties~\ref{thmstrongaH21}--\ref{thmstrongaH23} of Theorem~\ref{thmstrongaH2}. Finally, Property~\ref{thmstrongH24} is a direct consequence of Property~\ref{thmstrongaH24} of Theorem~\ref{thmstrongaH2} and the equation \eqref{nls}.
\medskip
\end{proof*}

\begin{thm}[\textbf{Finite time extinction and time decay estimates}]
\label{thmextH2}
Let Assumption~$\ref{ass1}$ be fulfilled with, in addition, $|\Omega|<\infty.$ Let $f\in W^{1,1}\big((0,\infty);L^2(\Omega)\big),$ $u_0\in H^1_0(\Omega)$ with $\Delta u_0\in L^2(\Omega)$ and let $u$ be the unique $H^2$-solution to \eqref{nls}--\eqref{u0} given by Theorem~$\ref{thmstrongH2}.$ Finally, assume there exists a finite time $T_0\ge0$ such that $f$ satisfies~\eqref{f}.
\begin{enumerate}
\item
\label{thmextH21}
If $N\le3$ then $u$ satisfies~\eqref{01} with,
\begin{gather}
\label{T*H2}
\frac{\|u(T_0\|_{L^2(\Omega)}^{1-m}}{(1-m)\Im(a)|\Omega|^\frac{1-m}2}+T_0\le
T_\star\le C\|u(T_0)\|_{L^2(\Omega)}^\frac{(1-m)(4-N)}4\|\Delta u\|_{L^\infty((0,\infty);L^2(\Omega))}^\frac{N(1-m)}4+T_0,
\end{gather}
for some $C=C(\Im(a),N,m).$
\item
\label{thmextH22}
If $N=4$ then for any $t\ge T_0,$
\begin{gather}
\label{thmrtdH21}
\|u(t)\|_{L^2(\Omega)}\le\|u(T_0)\|_{L^2(\Omega)}e^{-C(t-T_0)},
\end{gather}
where $C=C(\|\Delta u\|_{L^\infty((0,\infty);L^2(\Omega))},\Im(a),m).$
\item
\label{thmextH23}
If $N\ge5$ then for any $t\ge T_0,$
\begin{gather}
\label{thmrtdH22}
\|u(t)\|_{L^2(\Omega)}\le\dfrac{\|u(T_0)\|_{L^2(\Omega)}}{\left(1+C\|u(T_0)\|_{L^2(\Omega)}^\frac{(1-m)(N-4)}4(t-T_0)\right)^\frac4{(1-m)(N-4)}},
\end{gather}
where $C=C(\|\Delta u\|_{L^\infty((0,\infty);L^2(\Omega))},\Im(a),N,m).$
\item
\label{thmextH24}
If $N\le3$ then there exists $\eps_\star=\eps_\star(|a|,N,m)$ satisfying the following property. If
\begin{gather}
\label{thmextH2e1}
\begin{cases}
\|u_0\|_{L^2(\Omega)}^{2(1-\delta_2)}\le\eps_\star T_0,							\medskip \\
\|u_0\|_\star+\|f\|_{W^{1,1}((0,\infty);L^2(\Omega))}\le\eps_\star,				\medskip \\
\|f(t)\|_{L^2(\Omega)}^2\le\eps_\star\big(T_0-t\big)_+^\frac{2\delta_2-1}{1-\delta_2},
\end{cases}
\end{gather}
for almost every $t>0,$ where $\delta_2=\frac{m(4-N)+(4+N)}8\in\left(\frac12,1\right)$ and $\|u_0\|_\star^2=\|u_0\|_{H^1_0(\Omega)}^2+\|\Delta u_0\|_{L^2(\Omega)}^2,$ then $u$ satisfies~\eqref{01} with $T_\star=T_0.$
\end{enumerate}
\end{thm}

\begin{thm}[\textbf{Asymptotic behavior}]
\label{thmtdH2}
Let Assumption~$\ref{ass1}$ be fulfilled with $|\Omega|<\infty.$ Let $f\in W^{1,1}\big((0,\infty);L^2(\Omega)\big),$ $u_0\in H^1_0(\Omega)$ with $\Delta u_0\in L^2(\Omega)$ and let $u$ be the unique $H^2$-solution given by Theorem~$\ref{thmstrongH2}.$ Then,
\begin{gather}
\label{thmtdH21}
\vlim_{t\nearrow\infty}\|u(t)\|_{W^{1,q}(\Omega)}=\vlim_{t\nearrow\infty}\|u(t)\|_{L^p(\Omega)}=\vlim_{t\nearrow\infty}\frac{\d}{\d t}\|u(t)\|_{L^2(\Omega)}^2=0,
\end{gather}
for any $q\in(0,2]$ and $p\in\left(0,\frac{2N}{N-2}\right]$ $(p\in(0,\infty)$ if $N=2,$ $p\in(0,\infty]$ if $N=1).$
\end{thm}

\section{Proofs of the finite time extinction and asymptotic behavior theorems}
\label{proofext}

The proofs of Theorems~\ref{thmextN1} and \ref{thmextH2} are very close to those of the Theorems~3.5, 3.6, 3.7, 3.9, 3.11 and 3.12 in~\cite{MR4340780}. For convenience of the reader, we indicate the mains steps and refer to \cite{MR4340780} for more details.
\medskip \\
\begin{vproof}{of Theorems~\ref{thmextN1} and \ref{thmextH2}.}
By Gagliardo-Nirenberg's inequality, there exists $C_\GN=C(m,N)$ such that for any $v\in H^1_0(\Omega)\cap L^{m+1}(\Omega),$
\begin{align}
\label{GN1}
&	\|v\|_{L^2(\Omega)}^\frac{(N+2)-m(N-2)}2\le C_\GN\|v\|_{L^{m+1}(\Omega)}^{m+1}\|\nabla v\|_{L^2(\Omega)}^\frac{N(1-m)}2,	\\
\label{GN2}
&	\|v\|_{L^2(\Omega)}^\frac{(N+4)-m(N-4)}4\le C_\GN\|v\|_{L^{m+1}(\Omega)}^{m+1}\|\Delta v\|_{L^2(\Omega)}^\frac{N(1-m)}4,
	\text{ if also } \Delta v\in L^2(\Omega).
\end{align}
Now, suppose Assumptions~\ref{assN1} or the hypotheses of Theorem~\ref{thmextH2} are fulfilled. We choose $\ell=1$ for the proof of Theorems~\ref{thmextN1}, and $\ell=2$ for the proof of Theorems~\ref{thmextH2}. We let,
\begin{align*}
&	\delta_\ell=\frac{(N+2\ell)-m(N-2\ell)}{4\ell},		\quad	y(t)=\|u(t)\|_{L^2(\Omega)}^2, \; \forall t\ge0,	\\
&	\alpha=\Im(a)C_\GN^{-1},	\quad	\alpha_\ell=\alpha\|\nabla^\ell u\|_{L^\infty((0,\infty);L^2(\Omega))}^{-\frac{N(1-m)}{2\ell}},
						\quad	\nabla^2=\nabla.\nabla=\Delta.
\end{align*}
By \eqref{L2}, \eqref{GN1}--\eqref{GN2} and Hölder's inequality, we have for almost every $t\in(T_0,\infty),$
\begin{gather}
\label{edo}
y^\p(t)+2\alpha_\ell y(t)^{\delta_\ell}\le2\|f(t)\|_{L^2(\Omega)}y(t)^\frac12,	\\
y^\p(t)\ge-2\Im(a)|\Omega|^\frac{1-m}2y(t)^\frac{m+1}2.
\end{gather}
Using Assumptions~\ref{assN1} and the hypotheses of Theorem~\ref{thmextH2}, we obtain \eqref{01}--\eqref{thmrtdH12} and \eqref{T*H2}--\eqref{thmrtdH22} by integration (see also (2.10) in~\cite{MR4053613}). It remains to show the last property of the both theorems. By \eqref{thmstrongaH13}, there exists $\eps_\star=\eps_\star(|a|,m)$ with,
\begin{gather}
\label{eps*}
\eps_\star\le\min\left\{(2\delta_\ell-1)^{-\frac{2\delta_\ell-1}{\delta_\ell}}(\alpha\delta_\ell)^\frac1{1-\delta_\ell}(1-\delta_\ell)^\frac{2\delta_\ell-1}{\delta_\ell(1-\delta_\ell)},\alpha\,\delta_\ell\,(1-\delta_\ell)\right\},
\end{gather}
such that if \eqref{thmextN1e1} holds true then $\|\nabla u\|_{L^\infty((0,\infty);L^2(\Omega))}\le1.$ By \eqref{strongaH22}--\eqref{strongaH24}, \eqref{lemVL2} and \eqref{nls}, there exists $\eps_\star=\eps_\star(|a|,N,m)$ satisfying \eqref{eps*} such that under assumption \eqref{thmextH2e1}, we have $\|\Delta u\|_{L^\infty((0,\infty);L^2(\Omega))}\le1.$ Let $x_\star=(\alpha\delta_\ell(1-\delta_\ell)T_0)^\frac1{1-\delta_\ell}$ and $y_\star=\big(\alpha\delta_\ell^{\delta_\ell}(1-\delta_\ell)\big)^\frac1{1-\delta_\ell}.$ By \eqref{thmextN1e1}, \eqref{thmextH2e1} and \eqref{eps*},
\begin{gather}
\label{demthmextH2e1}
y(0)\le x_\star.
\end{gather}
Applying Young's inequality to \eqref{edo} and using \eqref{thmextN1e1}, \eqref{thmextH2e1} and \eqref{eps*}, we obtain
\begin{gather}
\label{demthmextH2e2}
y^\p(t)+\alpha y(t)^{\delta_\ell}\le y_\star\big(T_0-t\big)_+^\frac{\delta_\ell}{1-\delta_\ell},
\end{gather}
for almost every $t>0.$ By \eqref{demthmextH2e1}, \eqref{demthmextH2e2} and \cite[Lemma~5.2]{MR4053613}, $y(t)=0,$ for any $t\ge T_0.$
\medskip
\end{vproof}

\begin{vproof}{of Theorem~\ref{thm0w}.}
By \eqref{estthmweak} and density, we may assume that $f\in\Dr\big([0,\infty);L^2(\Omega)\big)$ and $u_0\in\Dr(\Omega).$ Then the result comes easily from Theorem~\ref{thmstrongaH2} and \eqref{L2}, by following the proof of \cite[Theorem~3.5]{MR4098330}.
\medskip
\end{vproof}

\begin{vproof}{of Theorem~\ref{thmtdH2}.}
Since $|\Omega|<\infty,$ we may assume that $q,p\ge2.$ Applying the proof of \cite[Theorem~3.14]{MR4340780}, the result follows.
\medskip
\end{vproof}

\section{Concluding remarks}
\label{scr}

\begin{enumerate}
\item
\label{scr1}
Do some $H^2$-solutions exist in the sense of~\cite{MR4098330,MR4053613,MR4340780} (see also \ref{rmkdefsol1} of Remark~\ref{rmkdefsol}) for $a\in D(m)$ $(0<m<1)$ but with $|\Omega|=\infty~?$
\item
\label{scr2}
In~\cite{MR4340780}, the existence of solutions is obtained with $m=0$ and $|\Omega|<\infty.$ The proof relies on the theory of maximal monotone operators on $L^2(\Omega)$ (Brezis~\cite{MR0348562}). Would it be possible to construct solutions but with $|\Omega|=\infty~?$ Of course, the method should be different since the nonlinearity $\frac{u}{|u|}$ does need not belong to $L^2(\Omega),$ and the notion of solutions might be revisited.
\item
\label{scr3}
The general method (that we shall call Method~1) to construct the solutions in \cite{MR4340780} is the following (in \cite{MR4098330}, the method is different and in \cite{MR4340780}, the domain $\Omega$ is bounded which makes the situation easier). We regularize the nonlinearity \eqref{gg} with \eqref{ge}. We associate operators $A$ and $A_\eps^m,$ to the nonlinearities \eqref{gg} and \eqref{ge}, respectively. We show that $(D(A_\eps^m),A_\eps^m)$ is maximal monotone in $L^2(\Omega).$ With help of a priori estimates, we may pass to the limit, as $\eps\searrow0,$ in the equation $(I+A_\eps^m)u_\eps=F$ to show that $(D(A),A)$ is maximal monotone in $L^2(\Omega).$ This permits to solve \eqref{nls} with initial data in $D(A),$ where, roughly speaking, $D(A)=H^2(\Omega)\cap H^1_0(\Omega)\cap L^{2m}(\Omega).$ The crucial tool to make such a choice of $D(A)$ possible is Lemma~4.2 in Bégout~\cite{MR4098330}. Another method which would be possible  (that we shall call Method~2) would be to show that $(D(A_\eps^m),A_\eps^m)$ is maximal monotone in $L^2(\Omega)$ and, with a priori estimates, to pass in the limit, as $\eps\searrow0,$ in the equation $\frac{\d u_\eps}{\d t}+A_\eps^mu_\eps=f(t,x),$ to solve \eqref{nls}. We then obtain the existence of $H^2$-solutions. With any of the two methods, the existence of $L^2$-solutions is obtained with help of a density argument and a result of continuous dependance such as Proposition~\ref{propdep}. Finally, $H^1_0$-solutions are obtained with a density argument and some a priori estimates obtained with help of \cite[Lemma~4.2]{MR4098330}. But when $a\in D(m),$ this lemma is no more valid. It follows that Method~1 fails to construct $H^2$-solutions, as well as Method~2 (actually, these both methods are equivalent). So we have to choose a larger domain $D(A)$ as \eqref{D}, which gives Theorem~\ref{thmstrongaH2} (by the way of Method~1), from which the existence of $L^2$-solutions follows. But due to the absence of a result such as in \cite[Lemma~4.2]{MR4098330}, we cannot establish estimates of the solution in the $H^1_0$-norm to construct $H^1_0$-solutions by density. This is why we apply Method~2 in this case. So, we may wonder if we might apply Method~2 from the beginning, without using Method~1. The answer is no because of the lack of a density result of smooth functions (roughly speaking, $H^2(\Omega)\cap H^1_0(\Omega))\cap L^{2m}(\Omega)$ is not dense in $D(A)$ defined by \eqref{D}). Finally, note that if we impose a stronger assumption of the initial data in Theorem~\ref{thmstrongaH2}, namely if we require that,
\begin{gather*}
u_0\in H^1_0(\Omega)\cap L^{2m}(\Omega) \text{ with } \Delta u_0\in L^2(\Omega),
\end{gather*}
instead of,
\begin{gather*}
u_0\in H^1_0(\Omega)\cap L^{m+1}(\Omega) \text{ with } \Delta u_0+a|u_0|^{-(1-m)}u_0\in L^2(\Omega),
\end{gather*}
then Method~2 completely works and we do not need to require to Method~1.
\end{enumerate}

\section*{Acknowledgements}
\baselineskip .5cm
P.~Bégout acknowledges funding from ANR under grant ANR-17-EURE-0010 (Investissements d’Avenir program). The research of J.~I.~D\'{\i}az was partially supported by the project ref. PID2020-112517GB-I00 of the DGISPI (Spain) and the Research Group MOMAT (Ref. 910480) of the UCM.

\baselineskip .4cm

\def\cprime{$^\prime$}

\end{document}